\documentclass[11pt]{article}
\usepackage{amssymb,amsmath,mathrsfs,amsfonts,graphicx,epsf} 
\usepackage{color}
\usepackage{xcolor}

 \input xy
 \xyoption{all}
 \newcommand{\ver}[1]{*+<2mm>[F-:<3pt>]{#1}}

 \usepackage{esint}

 \usepackage{fullpage}

\usepackage{authblk}

  \newcommand{\ffoot}[1]{}

\newcommand{\cut}{{\bf Cut}}

\newcommand{\bD}{\Delta}
\newcommand{\pphi}{\xi}

\newcommand{\ps}{\langle\cdot,\cdot\rangle}

\newcommand{\treint}{\sqint}
\newcommand{\f}{\mathfrak{f}}

\usepackage{mathrsfs}
\usepackage{graphicx}
\usepackage{amssymb}
\usepackage{color}
\newtheorem{theorem}{Theorem}
\newtheorem{corollary}{Corollary}
\newtheorem{lemma}{Lemma}
\newtheorem{proposition}{Proposition}
\newtheorem{definition}{Definition}
\newtheorem{remark}{Remark}

\newcommand{\bt}{\begin{theorem}}
\newcommand{\et}{\end{theorem}}
\newcommand{\bl}{\begin{lemma}}
\newcommand{\el}{\end{lemma}}
\newcommand{\bp}{\begin{proposition}}
\newcommand{\ep}{\end{proposition}}
\newcommand{\bc}{\begin{corollary}}
\newcommand{\ec}{\end{corollary}}
\newcommand{\bdeff}{\begin{definition}}
\newcommand{\edeff}{\end{definition}}
\newcommand{\brem}{\begin{remark}}
\newcommand{\erem}{\end{remark}}
\renewcommand{\r}[1]{(\ref{#1})}
\newcommand{\con}{{\mathcal C}}

\newcommand{\bi}{\begin{itemize}}
\newcommand{\iii}{\item}
\newcommand{\ei}{\end{itemize}}
\newcommand{\bd}{\begin{description}}
\newcommand{\ed}{\end{description}}

\newcommand{\bqn}{\begin{eqnarray}}
\newcommand{\eqn}{\end{eqnarray}}
\newcommand{\eqnn}{\nonumber\end{eqnarray}}

\newcommand{\nn}{\nonumber}
\newcommand{\ba}[1]{\begin{array}{#1}}
\newcommand{\ea}{\end{array}}

\newcommand{\R}{\mathbb{R}}

\newcommand{\lam}{\lambda}
\newcommand{\g}{\gamma}
\newcommand{\al}{\alpha}
\newcommand{\eps}{\varepsilon}

\newcommand{\om}{\omega}

\newcommand{\VecM}{\mathrm{Vec}(M)}

\newcommand{\Gq}{{\gg}_q}

\newcommand{\Zz}{\mathcal{Z}}
\renewcommand{\gg}{{\bf G}}
\newcommand{\sign}{\mathrm{sign}}
\newcommand{\E}{e}

\newcommand{\e}{\mbox{e}}

\newcommand{\HH}{{\bf (H0)}}

\newcommand{\frp}[2]{\frac{\partial #1}{\partial #2}}

\def\sign{\mathrm{sign}}

\title{\LARGE \bf
On almost-Riemannian surfaces
}

\author[*]{R.~Ghezzi} 
\affil[*]{Department of Mathematical Sciences, Rutgers University - Camden, 311 N 5th Street, Camden, NJ 08102, USA. {\tt roberta.ghezzi@rutgers.edu}}

\begin{document}

\maketitle

\begin{abstract}
An almost-Riemannian structure on a surface is a generalized Riemannian structure whose local orthonormal frames are given by Lie bracket generating pairs of vector fields that can become collinear. The distribution generated locally by  orthonormal frames has maximal rank at almost every point of the surface, but in general it has rank 1 on a nonempty set which is generically   a smooth curve. 
In this paper we provide a short introduction to 2-dimensional almost-Riemannian geometry highlighting its novelties with respect to Riemannian geometry. We present some results that investigate topological, metric and geometric aspects of almost-Riemannian surfaces from a local and global point of view.
\end{abstract}


\section{Introduction}

The purpose of this paper is to present a generalization of Riemannian geometry that naturally arises in the framework of control theory.
A Riemannian distance on a smooth surface  $M$ can be seen as the minimum-time function of 
an optimal control problem where admissible velocities are vectors of norm one. 
The control problem can be written locally as 
\bqn
\label{c-a}
\dot q=u X(q)+v Y(q)\,,~~~u^2+v^2\leq 1\,,
\eqn
where $\{X,Y\}$ is a local orthonormal frame.
Almost-Riemannian structures (ARSs for short) generalize Riemannian ones by allowing  $X$ and $Y$ to be collinear at some points.  In this case the corresponding Riemannian metric has singularities, but under generic conditions the distance is well-defined. For instance, if the two generators  satisfy the H\"ormander condition\footnote{The H\"ormander condition for a family ${\cal F}\subset\VecM$ of vector fields states that for every $q\in M$ Lie$_q{\cal F}=T_qM$.} (see for instance \cite{AgrBarBoscbook,bellaiche,jean1,jean2}), system \r{c-a} is completely controllable and the minimum-time function still defines a continuous distance  on the surface.

Let us denote by $\bD(q)$ the linear span of the two vector fields of a local orthonormal frame at a point $q$. Where $\bD(q)$ is 2-dimensional, the corresponding metric is Riemannian. Where $\bD(q)$ is 1-dimensional, the corresponding Riemannian metric is not well defined, but thanks to the H\"ormander condition one can still define the Carnot-Caratheodory distance between two points, which happens to be finite and continuous.
Generically, the singular set $\Zz=\{q\in M\mid \dim(\bD(q))=1\}$ is a 1-dimensional embedded submanifold  and there are three types of points: Riemannian points, Grushin points where $\bD(q)$ is 1-dimensional and dim($\bD(q)+ [\bD, \bD](q)) = 2$ and tangency points where dim($\bD(q) + [\bD, \bD](q)) = 1$ and the missing direction is obtained with one more Lie bracket. Generically, the following properties are satisfied: at Grushin points $\bD(q)$ is transversal to $\Zz$,  at tangency points $\bD(q)$ is tangent to $\Zz$ and tangency points are isolated.

Almost-Riemannian structures on surfaces were introduced in the context of hypoelliptic operators \cite{FL1,grushin1}. They appeared in problems of population transfer in quantum systems \cite{q1,q4,BCha} and they have applications to orbital transfer in space mechanics \cite{tannaka, BCa}.

The presence of a singular set enriches almost-Riemannian structures with several novelties with respect to the Riemannian case. The  aim of this paper is to present and discuss some of these aspects, mainly from a geometric point of view. 

For instance, as it happens in sub-Riemannian geoemetry, spheres centered at points of the singular set are never smooth and the asymptotic of the Carnot--Caratheodory distance is highly non-isotropic, see Section~\ref{sec:locres}. Moreover,  cut loci  accumulate at singular points, even in a non-smooth way, see \cite{BCGJ,kresta}. From the local point of view, the relations between curvature and conjugate points change, as the presence of a singular set permits the conjugate locus  to be nonempty even if the Gaussian curvature is negative, where it is defined. 
  From the global point of view, the relations between curvature and topology of the surface change as well.  Indeed,  Gauss--Bonnet-type formulas for ARSs were obtained \cite{ABS, high-order, euler}, where the role of the topology of the surface is instead played  by the topology of the components in which the singular set split the surface and by the contributions at tangency points, see Section~\ref{sec:gb}. As a consequence,   there exist  ARSs on surfaces  not necessarily parallelizable for which the integral of the curvature (suitably defined) vanishes, when in the Riemannian case the total curvature can vanish only on the torus. 
  Also, the location of the singular set and of tangency points plays a fundamental role when studying the Lipschitz equivalence class of the  Carnot--Caratheodory, see \cite{BCGS}. In particular, if in the Riemannian case all distances on the same surface are Lipschitz equivalent, in the almost-Riemannian case this is no longer true and the Lipschitz classification is finer than the differential one.
 Other interesting phenomena, not considered here, concerning the heat and the Schr\"odinger equation with the Laplace--Beltrami operator on an almost-Riemannian surfaces  were studied in  \cite{camillo}. In that paper it was proven that the singular set acts as a barrier for the heat flow and for a quantum particle, even if geodesics can pass through the singular set without singularities.

The structure of the paper is the following. In Section~\ref{sec:prel} we  recall the precise definition of almost-Riemannian structure on a surface and fix some notations. In Section~\ref{sec:grushin} we provide examples: the Grushin plane, for which we compute the optimal synthesis, and an example on the 2-sphere.
%
%
%
%
Section~\ref{loc} deals with local results. First in Section~\ref{locrep} we study local orthonormal frames  and state a result that classifies the types of points. Then, in Section~\ref{sec:locres} we focus on local properties  around tangency points such as the asymptotic of the distance as well as the cut and conjugate locus. Section~\ref{sec:forth} discusses functional invariants, i.e., functions that allow to recognize locally isometric structures.
In Section~\ref{glob} we address global problems. First we present in Section~\ref{sec:gb} several generalizations of the Gauss-Bonnet formula in the almost-Riemannian context. (For generalizations of Gauss--Bonnet formula in related contexts, see also \cite{agra-gauss, pelletier, pelle2}.) In Section~\ref{sec:topol} a characterization of ARSs using the topology of the vector bundle defining the structure is stated. Finally, in Section~\ref{sec:lipeq} we analyse  almost-Riemannian surfaces from a metric point of view presenting a classification result  with respect  to Lipschitz equivalence. 
We conclude  in Section~\ref{sec:riemvsariem} by some open questions.

\section{Basic definitions and notations}\label{sec:prel}

Let $M$ be a  smooth  connected surface without boundary and denote by $\VecM$ the  module of smooth vector fields on $M$. 
 Throughout the paper, unless specified,   objects are smooth, i.e.,  of class $\con^{\infty}$.

\bdeff\label{def:ars}
An \emph{almost-Riemannian structure} (ARS for short) on a surface $M$ is a triple $(E,\f,\ps)$, where (i) $E$ is a Euclidean  bundle   of rank two over $M$ (i.e. a vector bundle whose fibre is equipped with a smoothly-varying scalar product $\langle\cdot,\cdot\rangle_q$); (ii) $\f:E\rightarrow TM$ is  a morphism of vector bundles such that  $\f(E_{q})\subseteq T_{q}M$;  (iii) for every $q\in M$ Lie$_q\bD=T_qM$, where
$$
\bD=\{\f\circ\sigma\mid\sigma\mathrm{ \, section\,  of\, } E\}.
$$
 \edeff

ARSs on surfaces can be seen as a first generalization of Riemannian structures towards sub-Riemannian ones. Indeed, recall that a (constant rank) sub-Riemannian
structure on a manifold $N$ is given by $(D, g)$, where $D\subset TN$ is a sub-bundle of rank $k<\dim N$ such that Lie$_qD=T_qN$ for every $q\in N$ and $g$ is a Riemannian metric on $D$. When $\dim N=2$, the only possibility for a sub-bundle to be Lie bracket generating is that $D=TN$, i.e., there do not exist constant rank sub-Riemannian structures. Hence, in this case one needs to consider rank-varying modules $\bD\subset \mathrm{Vec} N$ instead of sub-bundles $D\subset TN$.

Let ${\mathcal S}=(E,\f,\langle\cdot,\cdot\rangle)$ be an ARS on a surface $M$.
The Lie bracket generating assumption in Definition~\ref{def:ars}, jointly with the fact that $\dim M=2$, implies that for each vector field $X\in\bD$ there exists a unique section $\sigma$ of $E$ such that $X=\f\circ\sigma$.  

Denote by $\bD(q)$  the linear subspace $\{V(q)\mid  V\in \bD\}=\f(E_q)\subseteq  T_q M$. 
 The set 
$$
\Zz=\{q\in M\mid \mathrm{dim }~\bD(q)<2\}
$$
is called \emph{singular set} and coincides with the set of points $q\in M$ at which $\f|_{E_q}$ is not injective, i.e.,  $\dim \bD(q)=1$. 

A property $(P)$ defined for 2-ARSs  is said to be {\it generic}
if for every rank-2 vector bundle $E$   over $M$, $(P)$ holds for every $\f$ in an  open and dense  subset of the set of  morphisms of vector bundles from $E$ to $TM$ inducing the identity on $M$, endowed with the 
$\con^\infty$-Whitney topology. 

The Euclidean structure $\ps$ on  $E$ induces a symmetric positive definite bilinear form $G(\cdot,\cdot)$ on the submodule $\bD$ as
$G(V,W)=\langle\sigma_V,\sigma_W\rangle$, where $\sigma_V, \sigma_W$ are the unique sections such that $V=\f\circ\sigma_V, W=\f\circ\sigma_W$.  At points $q$ where $\f|_{E_q}$ is an isomorphism $G$ acts as a tensor, i.e., $G(V,W)|_q$ depends only on $V(q),W(q)$. This is no longer true   at points belonging to   $\Zz$, which is generically  a smooth embedded submanifold of dimension 1. By definition, an ARS is Riemannian if and only if  $\f$ is an isomorphism of vector bundles or, equivalently, the singular set is empty.

If $(\sigma_1,\sigma_2)$ is an orthonormal frame for $\langle\cdot,\cdot\rangle$ on
 an open subset $U$ of $M$, an {\it  orthonormal frame for ${\cal S}$} on $U$ is 
 given by  $(\f\circ\sigma_1,\f\circ\sigma_2)$. 

For every $q\in M$  
and every $v\in\bD(q)$ define
$
\Gq(v)=\inf\{\langle u, u\rangle_q \mid u\in E_q,\f(u)=v\}$.
If $q\notin \Zz$, then $\Gq(v)=G(V,V)|_q$, for any vector field  $V$  such that $V(q)=v$. If $q\in\Zz$ then we have the inequality
$$
G(V,V)|_q\geq \Gq(V(q)).
$$

 An  absolutely continuous curve $\g:[0,T]\to M$ 
 is  {\it admissible} for ${\mathcal S}$ 
if  
there exists a measurable essentially bounded function 
$[0,T]\ni t\mapsto u(t)\in E_{\g(t)}
$ such that 
$\dot \g(t)=\f(u(t))$  for almost every $t\in[0,T]$. 
Given an admissible 
curve $\g:[0,T]\to M$, the {\it length of $\g$} is  
\bqn
\ell(\g)= \int_{0}^{T} \sqrt{ \gg_{\gamma(t)}(\dot \g(t))}~dt.
\eqnn
The {\it Carnot--Caratheodory distance} (or almost-Riemannian distance) on $M$  associated with 
${\mathcal S}$ is defined as
\bqn\nonumber
d(q_0,q_1)=\inf \{\ell(\g)\mid \g(0)=q_0,\g(T)=q_1, \g\ \mathrm{admissible}\}.
\eqn
The finiteness and the continuity of $d(\cdot,\cdot)$ with respect 
to the topology of $M$ are guaranteed by  the Lie bracket generating 
assumption  (see \cite[Theorem 5.2]{book2}).  
The Carnot--Caratheodory distance    endows $M$ with the 
structure of metric space compatible with the topology of $M$ as differential manifold.

Locally, the problem of finding a curve realizing the distance between two fixed points  $q_0,q_1\in M$ is naturally formulated as the distributional optimal control problem 
\bqn
\dot q= u_1 F_1(q)+u_2 F_2(q)\,,~~~u_i\in\R\,,
~~~\int_0^T 
\sqrt{u_1^2(t)+u_2^2(t)}~dt\to\min,~~q(0)=q_0,~~~q(T)=q_1,\eqnn
where $F_1,F_2$ is a local orthonormal frame for ${\cal S}$.

A {\it geodesic} for  ${\cal S}$  is an admissible 
curve $\g:[0,T]\to M$, such that 
 $\gg_{\gamma(t)}(\dot \g(t))$ is constant and   
for every sufficiently small interval 
$[t_1,t_2]\subset [0,T]$, $\g|_{[t_1,t_2]}$ is a minimizer of $\ell$. 
A geodesic for which $\gg_{\gamma(t)}(\dot \g(t))$ is (constantly) 
equal to one is said to be parameterized by arclength. 

Although  the metric tensor is singular at points of $\Zz$, geodesics are  well-defined and smooth. This can be proved by using classical methods of optimal control theory. 
The Pontryagin Maximum Principle \cite{pontryagin-book} provides a direct method to find geodesics for an ARS, as the H\"ormander condition ensures  the absence of abnormal extremals \cite[Proposition 1]{ABS}.  Indeed, it implies  that  an admissible curve parameterized by arclength is a geodesic if and only if it is the projection on $M$ of a solution of the Hamiltonian system corresponding to the Hamiltonian
\begin{equation}\label{eq:hamiltonian}
H(q,p)=\frac12( (p \cdot F_1(q))^2 + ( p \cdot F_2(q))^2),~~q\in U,\,p\in T^\ast_qU.
\end{equation}
lying on the level set  $H=1/2$, where $F_1,F_2$ is a local orthonormal frame for the structure on an open set $U$.  Notice that $H$ is well defined on the whole $T^*M$,  since  formula \r{eq:hamiltonian} does not depend on the choice of the orthonormal frame.
 When looking for a geodesic $\gamma$  realizing the distance from  a submanifold  $N$ (possibly of dimension zero), one should add the transversality condition $ p(0)\perp T_{\gamma(0)} N$.

The cut locus $\cut_N$ from  $N$ is the set of points $p$ 
for which there exists a geodesic realizing the distance between $N$ and $p$ that loses optimality after $p$.
 It is well known (see for instance \cite{agrompatto} for a proof in the three-dimensional contact case)
that, when there are no abnormal extremals, if $p\in \cut_{N}$ then one of the following two possibilities happen: (i) $p$ is reached optimally by more than one  geodesic; (ii) $p$ belongs to the first conjugate locus from $N$ defined as follows. To simplify the notation, assume that all geodesics are defined on $[0,\infty[$.
Define
\bqn
C^0=\{\lambda=(q,p)\in T^*M\mid q\in N,~H(q, p)=1/2,~p \perp T_{q}N\}
\eqnn
and 
\bqn
\exp:C^0\times[0,\infty[&\to& M\nn\\
(\lam,t)&\mapsto&\pi(e^{t \vec{H}}\lam)
\eqnn
where $\pi$ is the canonical projection $(q, p)\to q$ and $\vec{H}$ is the Hamiltonian vector field corresponding to $H$. The \emph{first conjugate time}  for the geodesic $\exp(\lambda,\cdot)$ is
\bqn
t_{conj}(\lam) = \min\{t > 0, (\lam,t) \mbox{ is a critical point of }\exp\}.
\eqnn
and the \emph{first conjugate locus from $N$} is $\{\exp(\lam,t_{conj}(\lam))\mid  \lam\in C^0 \}$.

\section{Examples}\label{sec:grushin}

\subsection{The Grushin plane}

To have a better understanding of the topic, let us study the simplest  
case of  genuinely 
almost-Riemannian surface.

 Consider  the ARS on $\R^2$ where $E=\R^2\times\R^2$, $\f((x,y),(a,b))=((x,y),(a, x b))$, and $\ps$ is the canonical Euclidean structure on $\R^2$.  In this case a global orthonormal frame is given by $F_1(x,y)=\partial_x$, $F_2(x,y)=x \partial_y$  and the 
singular set
is indeed nonempty, being equal to
the $y$-axis. 
This ARS is called \emph{Grushin plane},  named after  V.V.~Grushin who studied in \cite{grushin1} analytic properties  of  the operator $\partial^2_x+x^2\partial^2_y$ and of its multidimensional generalizations (see also \cite{FL1}).

The bilinear form $G(\cdot,\cdot)$ in coordinates $(x,y)$ reads
 $$
 \left(\ba{cc}
1 & 0\\ 0 & \frac 1 {x^2}
\ea\right),$$
and the Gaussian cuvature is $K(x, y) = -\frac{2}{x^2}$. Note that for every point outside the singular set the curvature is negative and it diverges to $-\infty$ as the point approaches $\Zz$.

The Grushin plane provides a very illustrative example as geodesics can be computed explicitly. Indeed, by the Pontryagin Maximum Principle, geodesics are projection on $\R^2$ of solutions of the Hamiltonian system associated with 
$$
H(x,y,p_x,p_y)=\frac1 2(p_x^2+x^2p_y^2).
$$
The system is
$$
\left\{\ba{ccc}
\dot x&=&p_x\\
\dot y&=&x^2 p_y\\
\dot p_x&=&-xp_y^2\\
\dot p_y&=&0
\ea\right.
$$
The last equation implies that along a geodesic $p_y(t)\equiv p_y(0)=a\in\R$.
 For $a=0$, 
  we obtain 
 $$
x(t)=x(0)+ p_x(0) t,~~y(t)\equiv 0, ~~~ p_x(t)\equiv p_x(0),
$$
that is, horizontal half-lines are geodesics.
For $a\neq 0$, the first and the third equation are the equations for an harmonic oscillator whose solution is
\begin{eqnarray}
x(t) &=& x(0)\cos(at) + \frac{\dot{x}(0)}{a}\sin(at)\label{eq:x}\\ p_x(t) &=& -ax(0)\sin(at) + \dot{x}(0)\cos(at).\nonumber
\end{eqnarray} 
Using the normalization $H=1/2$, the initial covector  $(p_x(0),a)$ satisfies
$$
H(x(0),y(0),p_x(0),p_y(0)) = \frac{1}{2} (p_x^2(0) + x^2(0)a^2) = \frac{1}{2},
$$
whence $$\dot{x}(0) = p_x(0) = \pm \sqrt{1 - x(0)^2a^2}.
$$ 
Integrating the equation for $y(\cdot)$ we get
\begin{eqnarray}
y(t) &=& y(0) + \frac{x(0)\dot{x}(0)}{2a} + \frac{a^2x(0)^2 + \dot{x}(0)^2}{2a} t + \nonumber\\&&+ \frac{a^2x(0)^2 - \dot{x}(0)^2}{2a}\sin(2at)  - \frac{x(0)\dot{x}(0)}{2a}\cos(2at).\label{eq:y}
\end{eqnarray}

\begin{figure}[h!]
\begin{center}
    \includegraphics[width=6cm]{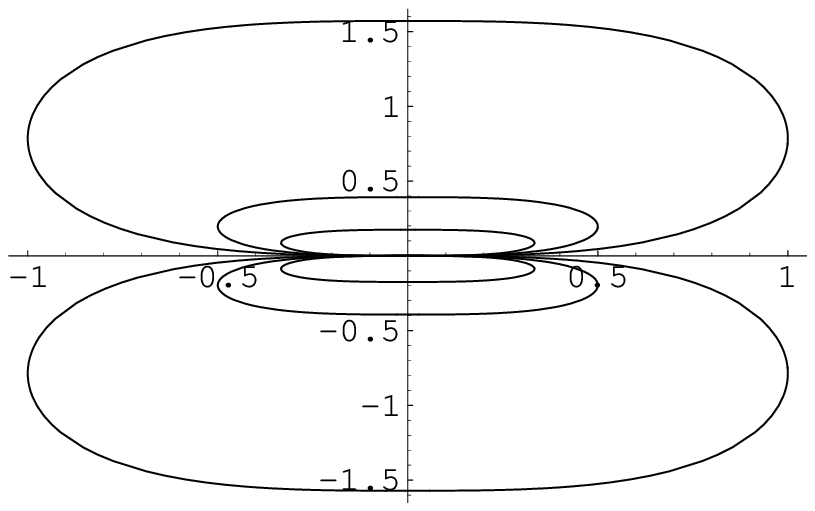}  
    \hspace{1.5cm}
    \includegraphics[width=6cm]{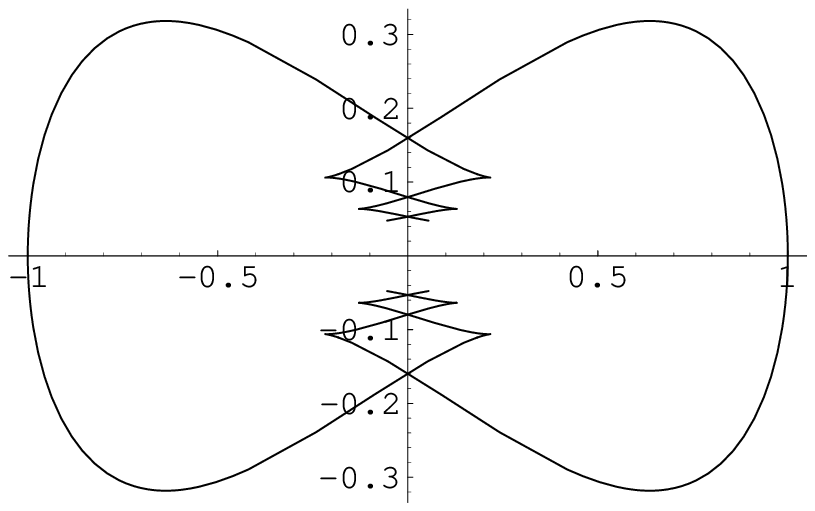} 
    
    a $~~~~~~\qquad\qquad\qquad\qquad~\qquad\qquad~~~~~~~~~~~~~~$ b

 \caption{Analysis of the Grushin plane at $(0,0)$. Figure~a shows some geodesics starting at $(0,0)$ corresponding to  $a = \pm 1,\pm 2,\pm 3$. Figure~b portraits the set of points whose distance from $(0,0)$ is equal to 1}\label{fronte}   
       \end{center}
    \end{figure}

Taking $(0,0)$ as starting point, the normalization condition $H=1/2~$ implies $~p_x(0)=1$ or $p_x(0)=-1$. When $p_x(0)=1$, respectively $p_x(0)=-1$, the geodesic enters the region $\{(x,y)\mid x>0\}$, respectively $\{(x,y)\mid x<0\}$.  
Hence, geodesics starting at $(0,0)$ 
have the form $(x^+(t,a),y(t,a))$ or $(x^-(t,a),y(t,a))$,   where
\begin{eqnarray}\label{geoori}
x^\pm(t,a) &=&\pm \frac{1}{a}\sin(at), \quad y(t,a) = \frac{1}{2a} t - \frac{1}{4a^2}\sin(2at) \quad a\neq 0,\\
x^\pm(t,0)&=& \pm t, \quad y(t,0)=0, \quad a=0.
\end{eqnarray}
Geodesics with $a>0$, respectively $a<0$, are contained in the half-plane $\{(x,y)\mid y>0\}$, respectively $\{(x,y)\mid y<0\}$. 
In Figure~\ref{fronte}a geodesics for some values of $a$ are portrayed, while Figure~\ref{fronte}b illustrates the set of points reached in time $t = 1$. Notice that this set is non-smooth. Also, the sphere\footnote{The sphere centered at $(x_0,y_0)$ of radius $r$ is defined as the boundary of the set $\{(x,y)\mid d((x,y),(x_0,y_0))<r\}$, where $d$ is the Carnot--Caratheodory distance. In general, for ARSs  this set contains $\{(x,y)\mid d((x,y),(x_0,y_0))=r\}$ as a proper subset.} centered at $(0,0)$ of radius $r=1$ is non-smooth. In contrast with what would happen in Riemannian geometry, this is the case for every positive radius, as it happens in constant-rank sub-Riemannian geometry. However, this is a consequence of the fact that the initial condition belongs to $\Zz$.
    
Let us compute the cut locus from $(0,0)$. The geodesics with $a=0$ 
 never lose optimality. Due to  the symmetries of the problem, 
 it is easy to see that the two geodesics $(x^+(t,a),y(t,a))$, $(x^-(t,a),y(t,a))$ are optimal until they intersect at time $t=\pi/|a|$. 
 As a consequence the cut locus from the origin is the set $\cut_{(0,0)}=\{(0,\alpha)\mid \alpha\in\R\setminus\{0\}\}$. Note that  $\cut_{(0,0)}$ accumulates at $(0,0)$. This is due to the fact that $(0,0)\in \Zz$ and represents another difference with the Riemannian case, where $\cut_p$ is always separated from $p$. 
 
\begin{figure}  
 \begin{center}
 \includegraphics[width=6cm]{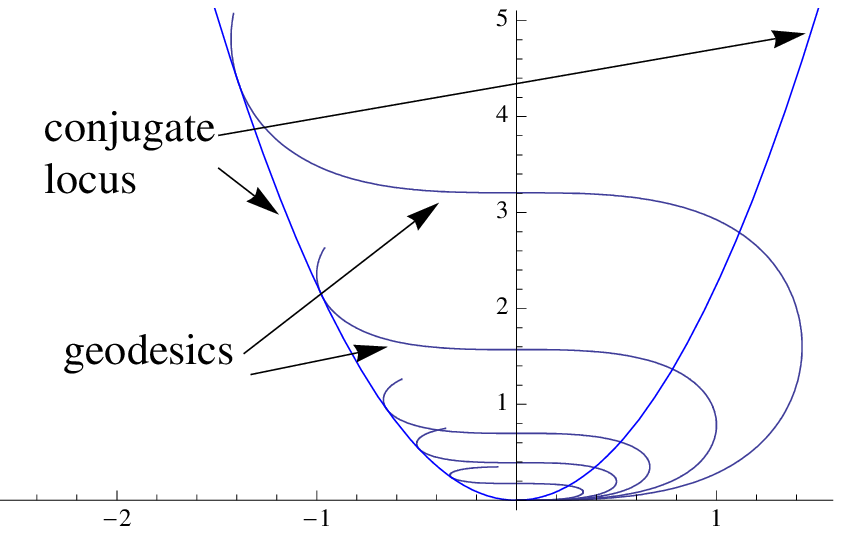} 
 \hspace{1.5cm}
 \includegraphics[width=6cm]{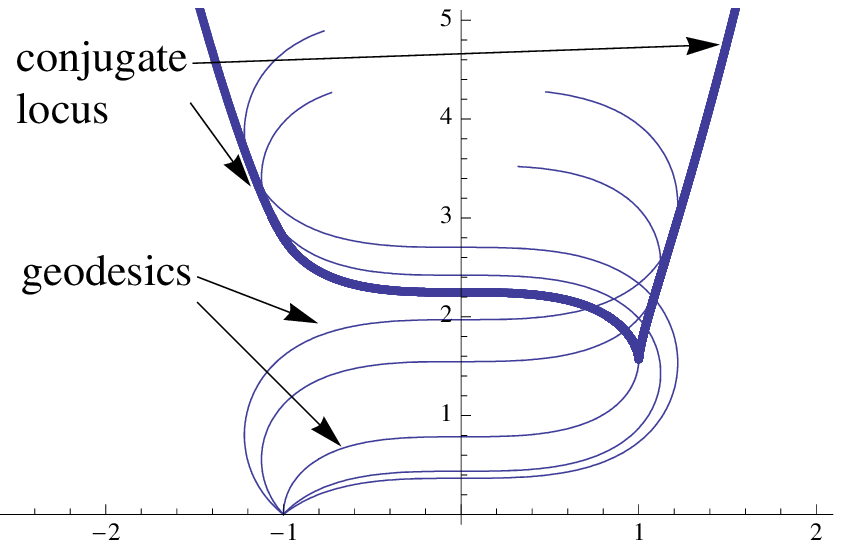} 
 
 a $~~~~~~\qquad\qquad\qquad~\qquad\qquad~~~~~~~~~~~~~~$ b
 \caption{Geodesics and upper part of first conjugate locus from a point in $\Zz$ (Figure~a) and from a point outside $\Zz$ (Figure~b)}\label{fig:conjbella}
   \end{center}
 \end{figure}

 To compute the conjugate locus from the origin, we find critical points of the map $\exp$ defined in
 Section~\ref{sec:prel}. Using \r{geoori} one can check that, for $a\neq 0$,  
 the first conjugate time on a geodesic $(x^\pm(t,a),y(t,a))$ coincides with the first positive root of the equation 
 $$
 a t \cos(a t)-\sin(a t)=0.
 $$ 
Therefore the conjugate locus from $(0,0)$ is the union of two parabolas 
$$
\left\{(x,y)\mid |y|=\frac{x^2}{2}\left(\frac{1}{\cos \tau \sin \tau}-\frac{1}{\tau}\right)\right\}\setminus\{(0,0)\},
$$
where $\tau$ is the first positive number such that $\tan \tau=\tau$. Note that since $(0,0)\in\Zz$, the conjugate locus accumulates at $(0,0)$.
 In particular, even if the curvature is negative for every point outside $\Zz$ the conjugate locus from $(0,0)$ is non-empty.  However, every geodesic reaches its conjugate time after crossing the singular set, see Figure~\ref{fig:conjbella}a.
 
\begin{figure}[t]
    \begin{center}
    \includegraphics[width=6cm]{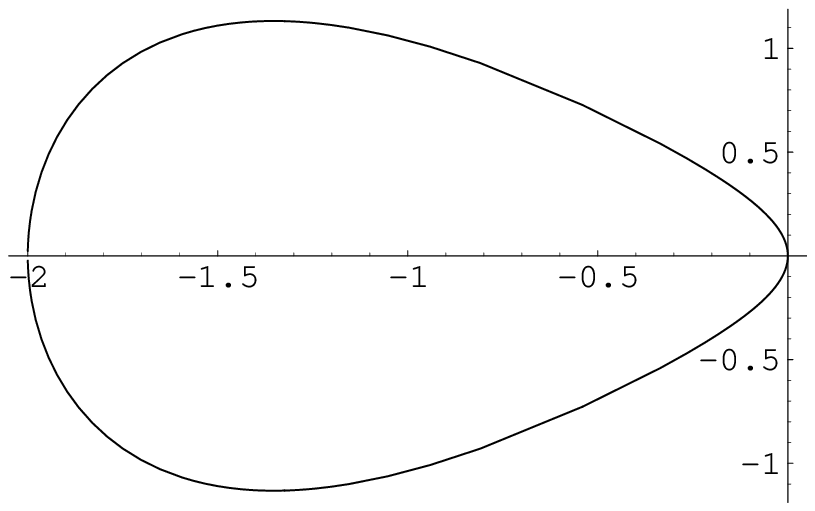}
    \hspace{1.5cm}
    \includegraphics[width=6cm]{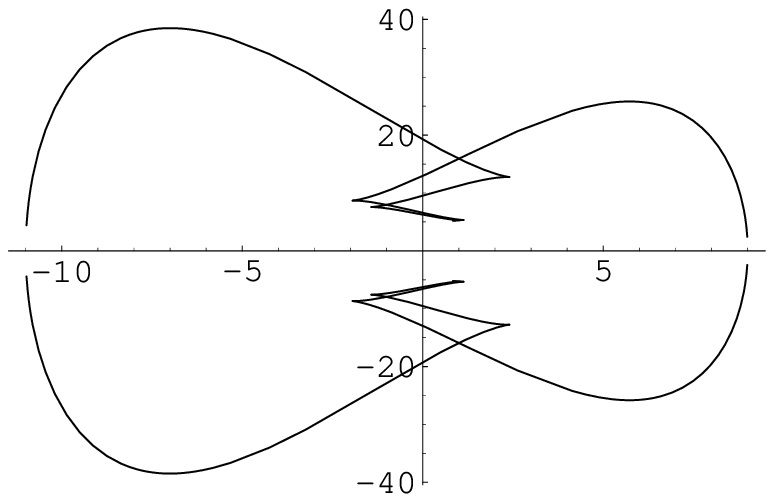}
\caption{The set of points $\{(x,y)\mid d((x,y),(-1,0))=r\}$, for $r=1, 10$.}\label{frontmenuno}\end{center}
    \end{figure}

  One may infer that the existence of conjugate points depends on  $(0,0)$ belonging to the singular set. However, this is not the case. Indeed, let us consider $(-1,0)$ as starting point. Let us parameterize the initial covector as $p_x(0)=\pm \sqrt{1- a^2}$,  $p_y(0)=a$. 
  When $a=0$ the solutions of the Hamiltonian system are $x^\pm(t,0)=-1\pm t$, $y(t,0)\equiv 0$. When $a\neq 0$,  using \r{eq:x}, \r{eq:y}, 
   the solutions are
   \begin{eqnarray*}
x^\pm(t,a) &=& \frac{-a\cos(at) \pm \sqrt{1 - a^2}\sin(at)}{a},\\
y^\pm(t,a) &=& \pm \frac{\sqrt{1 - a^2}}{a}(\cos(a t)^2-1) + \frac{2at - \sin(2at) + 2a^2\sin(2at)}{4a^2}.
\end{eqnarray*}
For $a\neq 0$, geodesics starting at $(-1,0)$ are curves of the form $(x^+,y^+), (x^+,-y^+)$ (when $p_x(0)=\sqrt{1-a^2}$) and of the form
$(x^-,y^-), (x^-,-y^-)$ (when $p_x(0)=-\sqrt{1- a^2})$. As one expects, since the metric is Riemannian at $(-1,0)$, for $r\leq 1$ the sphere centered at $(-1,0)$ of radius $r$ is smooth. This is no longer true for $r>d((-1,0),(0,0))=1$, see 
 Figure~\ref{frontmenuno}.
    
To compute the cut locus from $(-1,0)$, note first that the horizontal geodesics $(x^\pm(t,0),0)$ never lose optimality. On the other hand, when $a\neq 0$, one can check that the first time at which a geodesic
 $(x^+(t,a),y^+(t,a))$ intersects  another geodesic, namely $(x^-(t,a),y^-(t,a))$, is $t=\pi/|a|$. Similarly, the first intersection of $(x^+(t,a),-y^+(t,a))$ happens for $t=\pi/|a|$ with the geodesic $(x^-(t,a),-y^-(t,a))$. Hence, the cut locus from $(-1,0)$  is the union of two half-lines 
 $$
 \cut_{(-1,0)}=\{(1,\alpha)\mid\alpha\in(-\infty,-\pi/2)\cup (\pi/2,+\infty)\}
 .$$

As concerns the conjugate locus, computing critical points of $\exp$ it is easy to show that for each $a\neq 0$ every geodesic has a positive conjugate time. Hence, even if the curvature is negative and $(-1,0)\notin\Zz$, there exist conjugate points to $(-1,0)$, see Figure~\ref{fig:conjbella}b. Notice that each geodesic reaches the conjugate point after crossing at least one time the singular set.

\subsection{An example on the 2-sphere}\label{sec-sfera}

Consider another example on a compact surface. Let $M=S^2\subset\R^3$. Then every pair of vector fields on $S^2$  are linearly dependent on a nonempty set. Hence, if one wants to study a metric structure defined globally by a pair of vector fields on the 2-sphere, then one needs to consider an almost-Riemannian structure.

\begin{figure}
  \begin{center}
  \includegraphics[width=9cm]{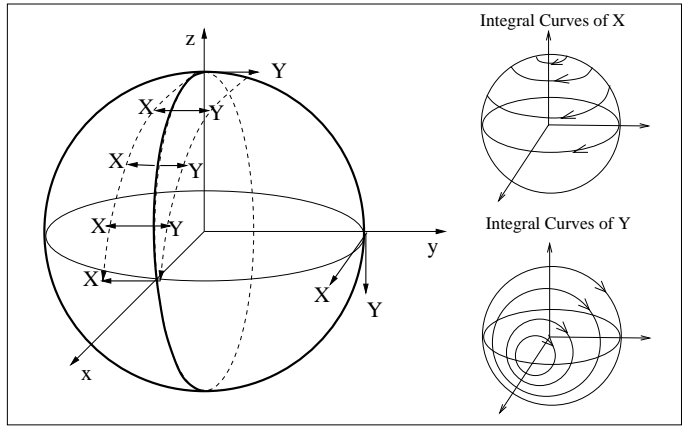}
  \caption{Almost-Riemannian structure on the 2-sphere}\label{f-sfera}
  \end{center}
  \end{figure}

Metric structures defined globally by a pair of vector fields on $S^2$  arise naturally in the context of quantum control (see \cite{q1,q4}). Indeed, consider the ARS on $S^2$ where $E$ is the trivial bundle of rank two over $S^2$ and the image under $\f$ of a global orthonormal frame for $\ps$ on $E$ is the pair  
$X(x,y,z)=(y,-x,0)$, $Y(x,y,z)=(0,z,-y)$. Then the two generators are linearly dependent on the  intersection of the sphere with the plane $\{(x,y,z)\mid y=0\}$  (see Figure~\ref{f-sfera}). 
In this model, the sphere represents a suitable state space reduction of a 
three-level quantum 
system and   the orthonormal generators $X$ and $Y$ are the infinitesimal rotations along two orthogonal axes,  
modeling 
the action on the system of
two lasers in the rotating wave approximation.  

Note that, as for the Grushin plane, the distribution is transversal to the singular set at each point. Indeed, this ARS is the compact correspondent to the Grushin plane.

\section{Local results}\label{loc}

\subsection{Local representations}\label{locrep}
The first important work  
studying general properties of ARSs  is \cite{ABS} where the authors provide  the characterization of generic ARSs by means of local 
representations.
\bdeff\label{def:locrep}
A  \emph{local representation} of an ARS ${\cal S}$ at  a point $q\in M$ is a pair of vector fields  $(X,Y)$ on $\R^2$ such that there exist: (i) a neighborhood $U$ of $q$ in $M$, a neighborhood $V$ of $(0,0)$ in $\R^2$ and a diffeomorphism $\varphi:U\rightarrow V$ such that $\varphi(q)=(0,0)$; (ii) a local orthonormal frame $(F_1,F_2)$ of ${\cal S}$ around $q$, such that  $\varphi_*F_1=X$, $\varphi_*F_2=Y$, where $\varphi_*$ denotes the push-forward.
\edeff
Under generic assumptions, it turns out that one can always construct a local representation $(X,Y)$ where the first vector field is rectified and the second one has a simple form.
The main assumption to obtain such result is the following. Set $\bD_1 = \bD$ and $\bD_{k+1}=\bD_k+[\bD,\bD_k]$, i.e., $\bD_k$ is the module spanned by Lie brackets of length less than $k$ between elements in $\bD$.
We say that $\cal S$ {\it satisfies condition} \HH\  if the following properties hold:
\bd
\iii[\HH] {(i)} $\Zz$ is an
embedded one-dimensional 
submanifold of
$M$;

\hspace{.1cm}{(ii)} the points $q\in M$ where $\bD_2(q)$ is
one-dimensional are isolated;

\hspace{.1cm}{(iii)}  $\bD_3(q)=T_qM$ for every $q\in M$.
\ed
A simple transversality argument allows to show that 
property \HH\ is generic for 2-ARSs. 
\brem\label{ho}
Throughout the paper, unless specified, we  always deal with ARSs satisfying \HH.
\erem
\begin{theorem}[\cite{ABS}]
\label{t-normal}
Given an ARS ${\mathcal S}$ on $M$, for every point
$q\in M$ there exist a local representation $(X,Y)$ of ${\cal S}$ at $q$ such that
$(X,Y)$
has one of the
forms
\bqn
\mathrm{(F1)}&& ~~X(x,y)=(1,0),~~~Y(x,y)=(0,e^{\phi(x,y)}), \nn  \\
\mathrm{(F2)}&& ~~X(x,y)=(1,0),~~~Y(x,y)=(0,x e^{\phi(x,y)}),\nn   \\
\mathrm{(F3)}&& ~~X(x,y)=(1,0),~~~Y(x,y)=(0,(y -x^2
\psi(x))e^{\pphi(x,y)}), \nn
\eqn
where $\phi$, $\pphi$ and $\psi$ are smooth real-valued functions such that
$\phi(0,y)=0$ and  $\psi(0)\neq 0$.
\et
 Thanks to Theorem~\ref{t-normal}, for a point $q\in M$ there are three possibilities. First, if $\bD(q)=T_qM$
then $q$ is said to be a
{\it Riemannian point} and  ${\mathcal
S}$ is locally described by (F1). 
Second, if $\bD(q)$ is one-dimensional and $\bD_2(q)=T_q M$ then we call  $q$ a  {\it Grushin
point} and
the local description (F2) applies. At Grushin points  $\bD(q)$ is transversal to $\Zz$ and the Lie bracket between the two elements of a local orthonormal frame is sufficient to span the tangent plane $T_qM$.
Third, if
$\Delta(q)=\Delta_2(q)$ has dimension one and $\bD_3(q)=T_q M$
then we say that $q$ is a {\it tangency point} and ${\mathcal
S}$ can be described near $q$ by  (F3). At tangency points the subspace $\bD(q)$ is tangent to $\Zz$ and the missing direction is obtained with a Lie bracket of length two between the two elements of a local orthonormal frame. We also set 
 $$
 {\cal T}=\{q\in \Zz\mid q \mbox{ tangency point of } {\cal S}\}.
 $$
Note that condition (ii) in assumption \HH\ ensures that tangency points are isolated, i.e., ${\cal T}$ is a discrete set.

\smallskip

The idea behind the proof of Theorem~\ref{t-normal} is that, since $\bD(q)$ is at least one-dimensional, given a local orthonormal frame for ${\cal S}$, one vector field can always be rectified. Then, to deduce the form of the other vector field one needs to construct a suitable coordinate system. To do this, the authors provide a procedure that allows to build a local representation starting from a smooth parameterized curve $c(\cdot)$ passing through the base point $q$ and transversal to the distribution at each point.
The transversality of $c(\cdot)$ to the distribution implies that the Carnot--Caratheodory distance from the support of $c(\cdot)$
 is  smooth on a neighborhood of $q$  (see \cite[Lemma 1]{ABS}). Given
a point $p$ near $q$, the first coordinate of $p$ is, by definition, 
 the distance between $p$ and the chosen curve,  with a suitable choice of sign. The second coordinate of $p$ is the 
 parameter corresponding to the point (on the chosen curve) that realizes the distance between $p$ and the curve (see Figure~\ref{fig-costr}). Using the inverse of the diffeomorphism defining this coordinate system, one can build a vector field belonging to $\bD$ of norm one and then complete it to a local orthonormal frame. The local orthonormal frame $(X,Y)$ obtained through this procedure has the form $X=(1,0)$, $Y=f(x,y) Y$, where $f$ is a smooth function.

\begin{figure}[h!]
\begin{center}
\input{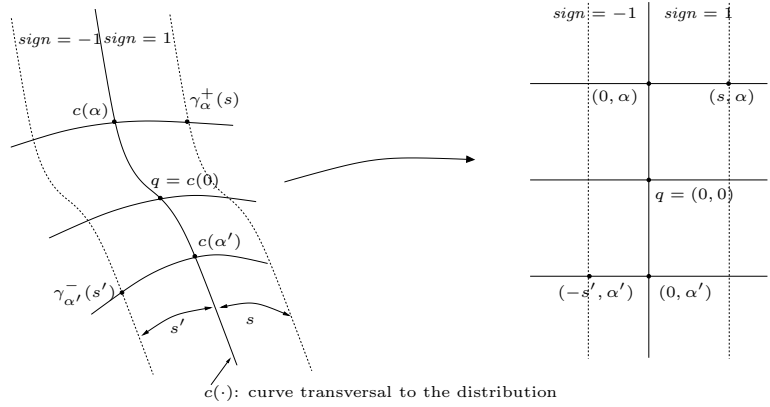}
\caption{The construction of coordinates starting from a parameterized curve $c:]-\varepsilon,\varepsilon[\to M$. We denote by $\g^{\pm}_\alpha$ the geodesic starting at $c(\alpha)$, parameterized by arclength, entering the region where {\it sign} $=\pm 1$ and such that $d(\g^{\pm}_\alpha(s),c(]-\varepsilon,\varepsilon[))=s$.  As the distribution is transversal to $c(\cdot)$ at each point, the distance from $c(]-\varepsilon,\varepsilon[)$ is smooth.}\label{fig-costr}
\end{center}
\end{figure}\subsection{Local analysis at tangency points} \label{sec:locres}

ARSs are characterized by the presence of a singular set, which includes Grushin and tangency points. These two types of points are essentially different from each other. 

At Grushin points, the local situation is described by the Grushin plane (see Section~\ref{sec:grushin}) which is the nilpotent approximation\footnote{The nilpotent approximation of a sub-Riemannian structure at a point is the analogue of the tangent space of a manifold. For the precise definition and analysis of the subject see for instance \cite{AgrBarBoscbook,bellaiche,jean2}.} of any ARS at a Grushin point. The fact that the distribution is   transversal to the singular set in a neighborhood of such points implies that the behaviour of the almost-Riemannian distance from $\Zz$ is similar on the two sides of $\Zz$, as we point out  in the discussion following Theorem~\ref{th:gb1} (see Section~\ref{sec:gb1}). 

At tangency points, the situation is more complicated  due to the fact that the asymptotic of the distance from the singular set is different from the two sides of the singular set. To see this, let us consider an example. Take the  ARS on $\R^2$   having  $X(x,y)=(1,0), ~Y(x,y)=(0,y-x^2)$ as a global orthonormal frame. In this case, the singular set is the parabola $\{(x,y)\mid y-x^2=0\}$ (see Figure~\ref{fig:astang}) and $(0,0)$ is a tangency point.
Define $M_\eps$ as the set of points $p\in\R^2$ such that $d(p,\Zz)>\eps$, where $d$ is the almost-Riemannian distance. Then $M_\eps$ is split in two connected components $M_\eps^>=M_\eps\cap\{(x,y)\mid y-x^2>0\}$ and  $M_\eps^<=M_\eps\cap\{(x,y)\mid y-x^2<0\}$ contained on the two different sides of $\Zz$. It turns out that for small values of $\eps$ both parts of the boundary of $M_\eps$ are non-smooth in a neighborhood of $(0,0)$ but have a corner. Also, the order at which they approach $\Zz$ as $\eps$ goes to zero is different. Indeed,  
taking another distance $\bar d$ on $\R^2$ associated with 
any Riemannian structure\footnote{For example choose $\bar d$ as the Euclidean distance on $\R^2$.},  one can check that
$
\bar d(M_\eps^>,\Zz)=O(\eps^2)
$, whereas $
\bar d(M_\eps^<,\Zz)=O(\eps^3).
$

\begin{figure}[h!]
\begin{center}
\includegraphics[width=6cm]{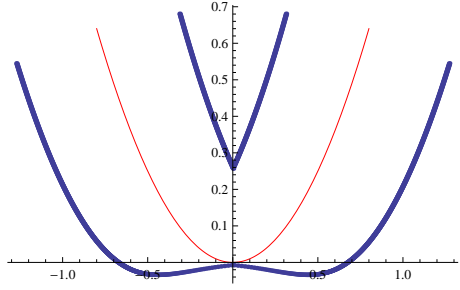}
\caption{The singular set (thick line) and the boundary of $M_\eps$ (bold line) in a neighborhood of a tangency point, for the structure on $\R^2$ whose orthonormal frame is $(1,0),(0,y-x^2)$ and $\eps=0.07$}\label{fig:astang}
\end{center}
\end{figure}
 
\smallskip

The asymmetric behaviour of the almost-Riemannian distance distinguishes tangency points from Grushin points and has various consequences. It affects the shape of the cut locus from the tangency point, as well as the shape of the cut locus from $\Zz$ near the tangency point. 
Moreover, together with the asymptotic of the Gaussian curvature, it is to be taken account of when proposing a notion of integrability with respect to an ARS (see Section~\ref{sec:gbtang}).  
Let us analyse singular loci  around tangency points.

%
%
Consider  an ARS in a neighborhood of a tangency point. Thanks to Theorem~\ref{t-normal}, this is equivalent to taking the  ARS  on $\R^2$ whose orthonormal frame is 
\bqn\label{eq-f3}
X(x,y)=(1,0),~~~Y(x,y)=(0,(y -x^2
\psi(x))e^{\pphi(x,y)}),
\eqn
where $\psi,\pphi$ are smooth functions depending on the structure and  $\psi(0)\neq 0$. 
 It is easy to see that     the nilpotent approximation of \r{eq-f3} at $(0,0)$  is the ARS defined by the orthonormal frame
\bqn\label{nilp-f3}
\hat X(x,y)=(1,0),~~~\hat Y(x,y)=(0,-\gamma x^2),
\eqn
where  $\gamma=\psi(0)e^{\xi(0,0)}$. The singular set of this structure is the $y$-axis and at each singular point $\dim\hat\bD(0,y)=1$. Therefore,  this structure does not satisfy condition (ii) in \HH.
The optimal synthesis for this ARS was computed explicitly in \cite{ag-bonn, BC}  in terms of Jacobi elliptic functions. Even though such synthesis does not reflect qualitative properties of the one for the generic case \r{eq-f3}, it can be used to compute the jet of the exponential map for the ARS \r{eq-f3}.
%
%
%
%
%
The development at $(0,0)$ of the orthonormal frame in \r{eq-f3} truncated at order zero\footnote{The coordinate functions $(x,y)$ have weights $(1,3)$, see \cite{bellaiche}.} is
 \begin{equation}\label{AR0}
\tilde{X}=\frp{}{x}, \tilde{Y}= \gamma(y-x^2-\eps' x^3)\frp{}{y},
\end{equation}
where $\eps'=\psi'(0)+\psi(0)\xi_x(0,0)$.
The following proposition computes the exponential map at the origin for the ARS defined in \r{AR0}. It happens that  higher order terms in the expansion of the elements of the orthonormal frame in \r{eq-f3} do not affect the estimation of  the exponential map and, consequently,  the order zero is sufficient to describe the cut and conjugate loci from the tangency point, at least qualitatively. 

\bp[\cite{BCGJ}] \label{proposition31}
Consider the ARS on $\R^2$ defined by the orthonormal frame given in \r{AR0}. The solution of the Hamiltonian system associated with 
$$
H(x,y,p_x,p_y)=\frac{1}{2}(p_x^2+\gamma^2(y -x^2-\eps' x^3)^2p_y^2)
$$ with initial condition 
$(x,y,p_x,p_y)|_{t=0}=(0,0,\pm 1,a)$
with $|a|\sim +\infty$ can be expanded as
\begin{eqnarray*}
x(t,\eta) & = & \eta~ x^{0}(t/\eta)+\eta^2x^{1}(t/\eta)+o(\eta^2),\\
y(t,\eta) & = & \eta^3~y^{0}(t/\eta)+\eta^4 y^{1}(t/\eta)+o(\eta^4),
\end{eqnarray*}
where $\eta=\frac{1}{\sqrt{|a|}}$, $(x^0,y^0,p_x^0,p_y^0)$ is the extremal of the nilpotent approximation\footnote{that is, $(x^0,y^0,p_x^0,p_y^0)$ is a solution of the Hamiltonian system associated with $H^0=1/2(p_x^2+x^4p_y^2)$, see \r{nilp-f3}} with initial condition $(x^0,y^0,p_x^0,p_y^0)|_{t=0}=(0,0,\pm 1,\sign(a))$, and 
\bqn
\left\{
\ba{lll}
\dot{x}^1&=&p_{x}^{1},\\
\dot{y}^1&=&\gamma^2(p_y^1(x^{0})^4+4 p_{y}^{0}(x^0)^3x^1-2p_y^0( (x^{0})^{2}y^0-\eps'(x^0)^5)),\\
\dot{p}^1_x&=& -\gamma^2(4 p_y^0 p_y^1 (x^0)^3+6(p_y^0)^2(x^0)^2x^1-(p_y^0)^2(2 x^0 y^0-5\eps' (x^0)^4))\\
\dot{p}^1_y&=& \gamma^2{p_{y}^{0}}^{2}{x^{0}}^{2},
\ea\right.
\eqnn
with initial condition $(x^1,y^1,p_x^1,p_y^1)|_{t=0}=(0,0,0,0)$.
\ep

Proposition~\ref{proposition31} was used to estimate the  conjugate locus  from $(0,0)$ for the  ARS on $\R^2$ defined by the orthonormal frame \r{eq-f3}, see \cite[Proposition 5]{BCGJ}. There exists  a constant  $\alpha\neq 0$ such that the conjugate locus from $(0,0)$  accumulates at $(0,0)$ as the set  
$$
\{(x,y)\mid  y=\alpha x^{3}\}\cup \{(x,y)\mid   y=-\alpha x^{3}\}\setminus\{(0,0)\}.
$$
 
 The shape of the cut locus from a tangency point (see Figure~\ref{figurone}) is described by the following result.
\bp[\cite{BCGJ}]\label{formacut}
Let ${\cal S}$ be an ARS on $M$ and $q\in M$ be a tangency point  such that there exists a local representation of the type (F3) for ${\cal S}$ at $q$ with the property
$$
\psi'(0)+\psi(0)\partial_x\xi(0,0)\neq0.
$$
 Then the cut locus from the tangency point  accumulates at   $q$ as an asymmetric cusp whose branches are  separated locally by $\Zz$.
In the coordinate system where the chosen local representation  is (F3), the cut locus accumulates as the set
$$
\{(x,y)\mid y>0,  y^{2}-\alpha_1 x^{3}=0\}\cup \{(x,y)\mid y<0,  y^{2}-\alpha_2 x^{3}=0\},
$$
with
$\al_i=c_ie^{2\xi(0,0)}/(\psi'(0)+\psi(0)\partial_x\xi(0,0))^3$, the constants $c_i$ being nonzero and independent on the structure.
\ep
Note that both the conjugate locus and $\cut_{(0,0)}$  accumulates at $(0,0)$ with tangent direction parallel to the distribution. This is no longer true for $\cut_\Zz$  in a neighborhood of a tangency point.  
A description of such locus is given by the following theorem, see also Figure~\ref{figurone}.

\bt[\cite{kresta}]\label{prop-cutz}
Let ${\cal S}$ be an ARS on $M$ and $q\in M$ be a tangency point  such that there exists a local representation of the type (F3) for ${\cal S}$ at $q$ with the property
$$
\eps'=\psi'(0)+\psi(0)\partial_x\xi(0,0)\neq0.
$$
Then the cut locus from the singular set $\Zz$  accumulates at   $q$ as the union of two curves locally separated  by $\Zz$. 
 One of them is contained in the set $\{y>x^2\psi(x)\}$,   takes the form
$$
\{(-1/2\psi'(0)t+o(t),t+o(t))\mid t>0\},
$$
and accumulates at $q$ transversally to the distribution.
The other one is contained in the set $\{y<x^2\psi(x)\}$ and takes the form
$$
\{(\eps' \, \om \, t^2+o(t^2),- t^3+o(t^3))\mid t>0\},
$$
where $\om\neq 0 $ is a constant depending on the structure. This part of the cut locus accumulates at $q$ 
with
 tangent direction at $q$ belonging to the distribution. 
%
%
\et
\begin{figure}[h!]
\begin{center}
\includegraphics[width=6cm]{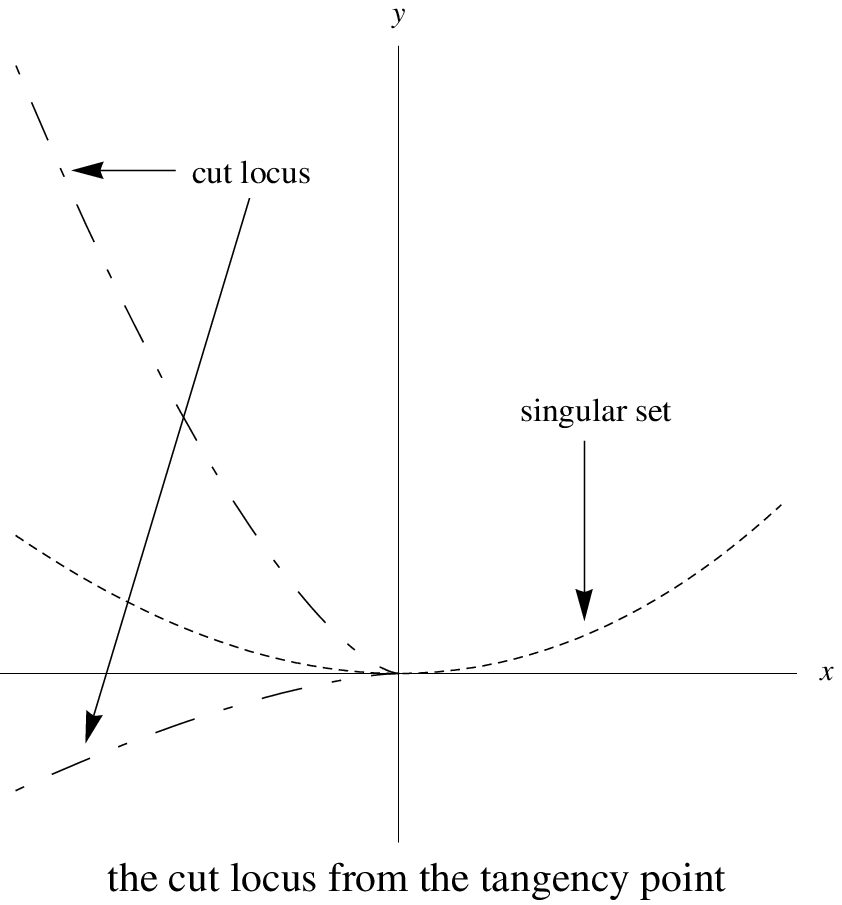}
\hspace{2cm}
\includegraphics[width=6cm]{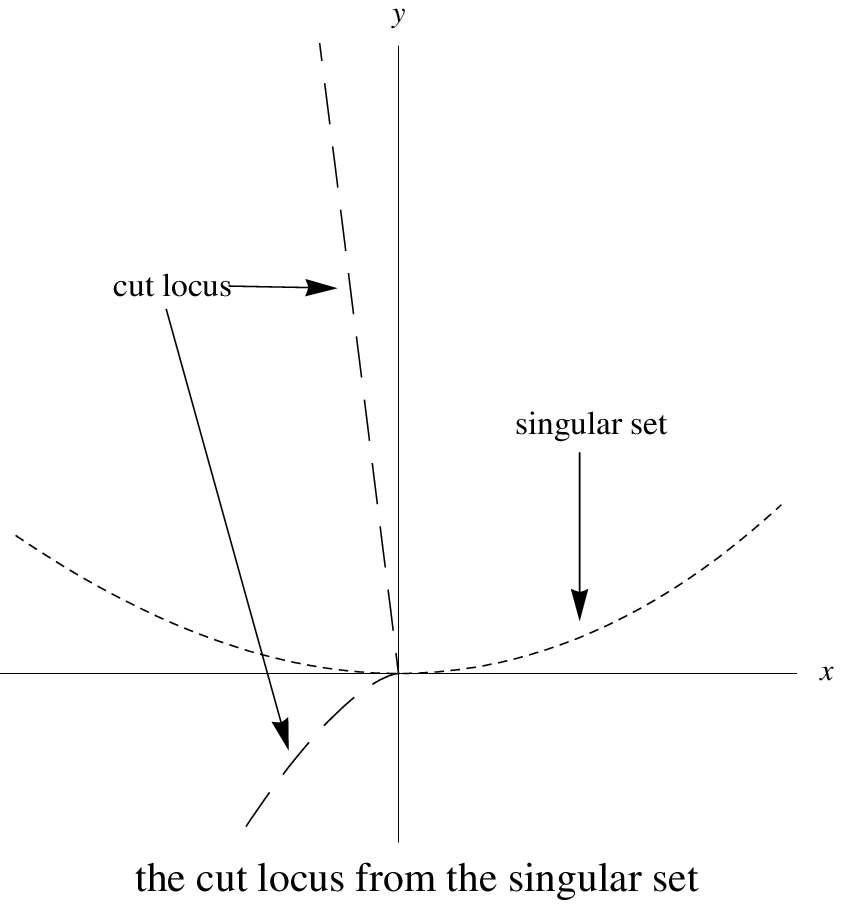}

\vspace{.5cm}

\includegraphics[width=6cm]{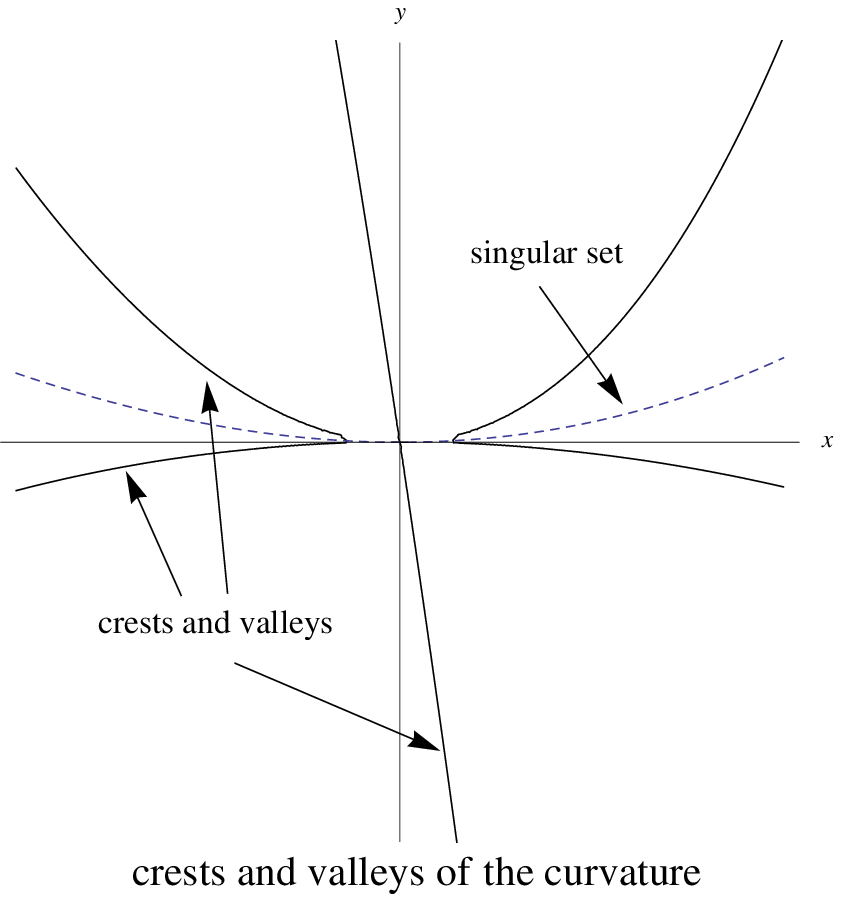}
\hspace{2cm}
\includegraphics[width=6cm]{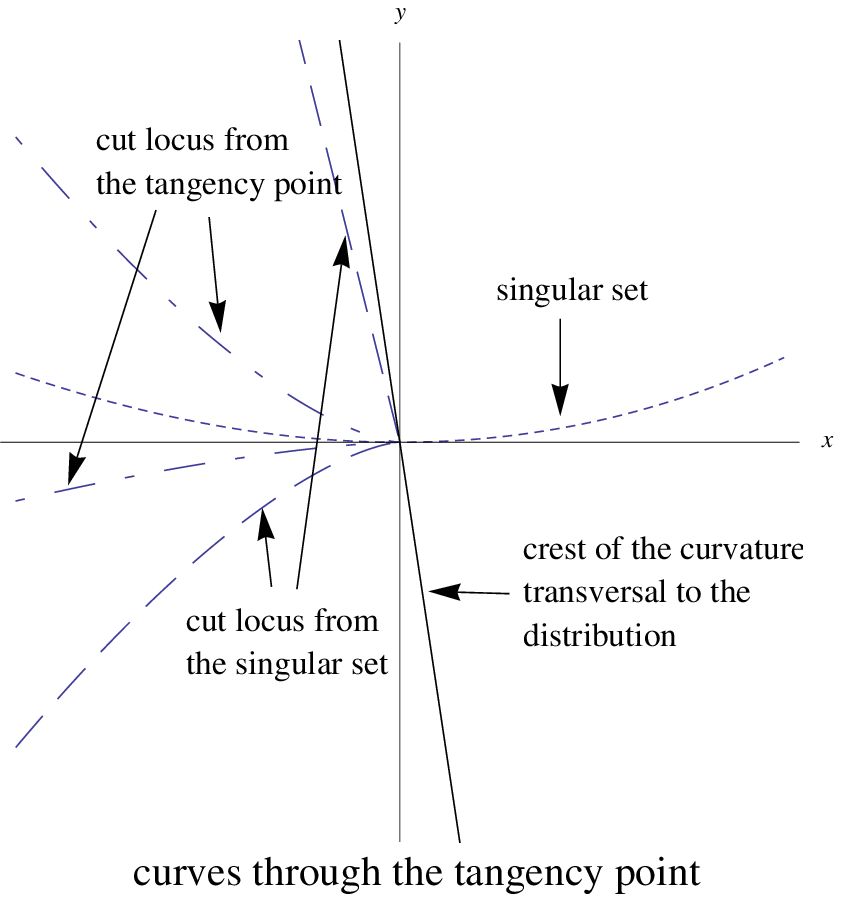}
\end{center}
\caption{The singular set (dotted line), the cut locus from a tangency point 
(semidashed line),
the cut locus from the singular set (dashed line), and the set   of crests and valleys of $K$  (solid lines)  
for the ARS with orthonormal frame $F_1=\frp{}{x},\, F_2=(y-x^2-x^3)\frp{}
{y}$. In this case  there are three curves in the set of crests and valleys of the curvature,  only one of which is transversal to  the distribution.}
\label{figurone}
\end{figure}

%
%
%

\subsection{Functional invariants}\label{sec:forth}

  In a recent paper \cite{kresta} the authors address the problem of  finding normal forms for ARSs that are completely reduced, in the sense that they depends only on the  ARS and not on its local representation.
This consists in finding a canonical choice
for a local coordinate system  and for a local orthonormal frame, i.e., two vector fields $X$ and $Y$ defined in a neighborhood of the origin on $\R^2$. 
Once 
a canonical choice of $(X,Y)$ is provided, the two components of  $X$ and the two components of $Y$ are functional invariants of the structure, in the sense that locally isometric structures have the same components. Moreover, they   permits to recognize locally isometric structures: if two structures have the same invariants in a neighborhood of a point, then they are locally isometric.
Notice that the problem of finding a set of invariants that determines the structure up to local isometries is not completely trivial even in the simplest case of Riemannian points. See for instance the discussion in \cite{zelclass, kulkarni}.  Indeed even if one is able to fix canonically a coordinate system, the Gaussian curvature in that coordinate system is an invariant, but  there are non-locally isometric structures having the same curvature.

A first step in finding normal forms is Theorem~\ref{t-normal}. However, the local representations corresponding to Riemannian and tangency points  are not completely reduced. 
Indeed, there  exist changes of coordinates and rotations  of the frame for which an orthonormal basis has the same expression as in (F1)  (respectively (F3)), but with a different function $\phi$ (respectively with different functions $\psi$ and $\xi$). 
Recall that to build the local representations in Theorem~\ref{t-normal} a specific procedure based on the choice of a smooth parameterized curve transversal to the distribution was used.
 If the parameterized curve used in this  construction can be canonically built, then one gets 
 a local representation of the form $X=(1,0)$, $Y=(0,f)$ which  cannot be further reduced. Hence,  $f$ is a functional invariant that completely determines the structure up to local isometries.   In  \cite{kresta}, a canonical choice of the parameterized curve at each point is provided and the properties of the functional invariant $f$ are studied.

For Riemannian points, a canonical parameterized curve transversal to the distribution can be easily identified, at least at  points where the gradient of the Gaussian curvature is non-zero: one can use the level set of the curvature passing from the point, parameterized by arclength. For points where the gradient of the curvature vanishes, under additional generic conditions,   a smooth parameterized canonical curve passing through the point can be selected among the crests and valleys of the curvature, see \cite[Sections 4.1.3, 4.1.4]{kresta}.
 
  For Grushin points, a canonical curve transversal to the distribution is the set $\Zz$. This curve has also a natural parameterization  and was used in \cite{ABS} to get the local representation (F2) that, as a consequence, cannot be further reduced.
  
To obtain the  local expression (F3),  
the choice of the smooth parameterized curve was arbitrary and not canonical. 
 The analysis in \cite{kresta}  is aimed at finding a canonical one and, as a consequence, obtain a   functional invariant at a tangency point that completely determines the structure. 
The most natural candidate for such a curve is the cut locus from the tangency point. Neverthless,
by Proposition~\ref{formacut} this is not a good choice, as  in general   
the cut locus  from the point is not smooth but has an asymmetric cusp (see Figure~\ref{figurone}). 
Another possible candidate is the cut locus from the singular set in a neighborhood 
of the tangency point, but  Theorem~\ref{prop-cutz} states that
the cut locus from $\Zz$ is non-smooth in a neighborhood of a tangency point, see Figure~\ref{figurone}. 
A third possibility is to look for curves which are crests or
valleys of the Gaussian curvature and  intersect transversally the singular set at a tangency point. The main result in \cite{kresta}   consists in the  proof of  the existence of  such a curve.  Moreover, this curve admits a canonical regular  parameterization.
%
%

 An example of a crest of the curvature at a tangency point is shown in Figure~\ref{fig-KRESTACCIA}.
 \begin{figure}[h!]
\begin{center}
\includegraphics[width=8cm]{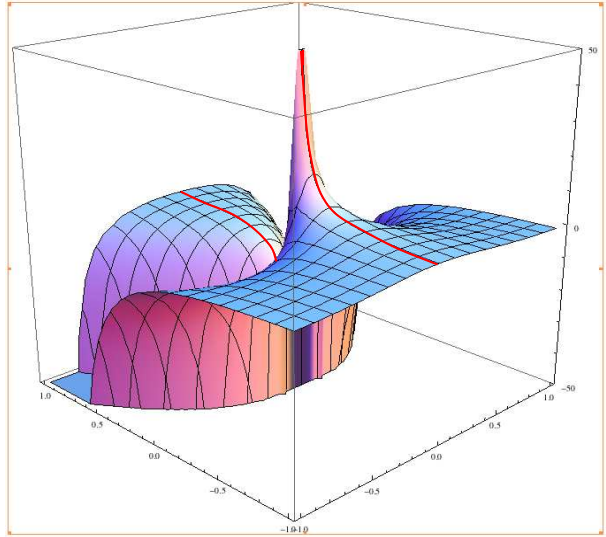}
\caption{The graph of the Gaussian curvature $K=\frac{-2(3 x^2 +y)}{(y-x^2)^2}$ for the almost-Riemannian structure on $\R^2$ having 
$X(x,y)=(1,0),\, Y(x,y)=(0, y-x^2)$ as orthonormal frame. The curvature has a crest  passing through the tangency point $(0,0)$}\label{fig-KRESTACCIA}
\end{center}
\end{figure}

\section{Global results}\label{glob}
 
In this section we address global questions for almost-Riemannian surfaces.
The general context in which all the results can be stated is the one of totally oriented structures.

\bdeff
\label{d-oriented}
An ARS is said to be {\em oriented} if $E$ is oriented as vector bundle. We say that an ARS is {\em totally oriented} if both $E$ and $M$ are oriented. For a totally oriented ARS, $M$ is split into two open sets $M^+$, $M^-$ such that $\Zz=\partial M^+=\partial M^-$, $\f:E|_{M^+}\rightarrow TM^+$ is an orientation-preserving isomorphism and $\f:E|_{M^-}\rightarrow TM^-$ is an orientation reversing-isomorphism.
\edeff

\brem\label{assunzioni}
Throughout this section we deal with totally oriented ARS on compact surfaces (satisfying assumption \HH, see Remark~\ref{ho}). 
\erem
\subsection{Gauss-Bonnet formulas}\label{sec:gb}
 
 The first global result for almost-Riemannian surfaces is a Gauss--Bonnet-type formula proved in \cite{ABS}. Such theorem gives a first sight of how the relation between curvature and topology changes in the almost-Riemannian context. Indeed, not only the topology of the surface shows up in the Gauss-Bonnet formula, but also the way the singular set embeds in the surface. Moreover, a central role is played by  the topology of the vector bundle, see Corollary~\ref{th:totale}.
 
 Other generalizations of the Gauss--Bonnet formula can be found in \cite{agra-gauss} for contact sub-Riemannian manifolds and in \cite{pelletier,pelle2} for pseudo-Riemannian manifolds.

\subsubsection{ARSs without tangency points}\label{sec:gb1}
 
The Gauss--Bonnet formula for a compact oriented Riemannian surface $M$ states that
 \begin{equation}\label{eq:gbr}
 \int_MK dA=2\pi\chi(M),
 \end{equation}
 where $K$ is the Gaussian curvature and $dA$ is the volume form associated with the metric.

 The first step to generalize \r{eq:gbr} to ARSs is to introduce an analogous of the volume form $dA$.  
 To this aim, the idea is to take a non-degenerate volume form $\omega$ on $E$ for the Euclidean structure and push it forward on $M$ using the morphism $\f$. Since $\f$ is not an isomorphism on $\Zz$, this does not give rise to a two form on the whole surface $M$ but only on the set $M\setminus\Zz$. Let  $dA_s=f_*w$ denote such a form. Then  $dA_s$ defines on $M^+$ the same orientation as the one chosen on $M$, the opposite one  on $M^-$.

By construction, $dA_s\in \Omega^2(M\setminus \Zz)$ and it diverges when approaching the singular set. Also, the Gaussian curvature diverges when approaching to $\Zz$. 
To overcome these issues and provide a way of integrating $K dA_s$ over $M$ the idea in \cite{ABS} is to integrate $K dA_s$ over the set of points $\eps$-far from the singular set and then let $\eps$ go to zero. More precisely, let
\begin{equation}\label{eq:meps}
M_\eps=\{q\in M\mid d(q,\Zz)>\eps\},
\end{equation}
 (where $d$ is the almost-Riemannian distance). Then the Gaussian curvature is said to be \emph{integrable with respect to the ARS} if the limit 
 \begin{equation}\label{eq:lim}
 \lim_{\eps \to 0}\int_{M_\eps}K dA_s 
 \end{equation}
 exists. In this case such limit is denoted by $\int_M K dA_s$.
 
 It turns out that if there are not tangency points then the limit in \r{eq:lim} exists and can be calculated in terms of the topology of the ARS.
 
 \bt[\cite{ABS}]\label{th:gb1}
If there are no tangency points, then 
 \begin{equation}\label{eq:gb1}
 \int_M K dA_s=2\pi (\chi(M^+)-\chi(M^-)).
 \end{equation}
 \et
 
 To explain the assumption of absence of tangency points, let us go through the main ideas in the proof of \r{eq:gb1}. By definition,
 $$
 \int_{M_\eps} K dA_s=\int_{M_\eps\cap M^+} K dA_s+\int_{M_\eps\cap M^-}K dA_s.
$$ 
 On the sets $M_\eps\cap M^{\pm}$ we can apply the classical Gauss--Bonnet formula to get
 $$
 \int_{M_{\eps}\cap M^+}K dA_s=2\pi\chi(M_\eps\cap M^+)-\int_{\partial(M_\eps\cap M^+)}k_g d\sigma,
 $$
 where $k_g$ is the geodesic curvature and $\partial(M^\eps\cap M^+)$ carries the orientation induced by $M^\eps\cap M^+$. Taking account of orientations, similarly we obtain 
 $$
 \int_{M_{\eps}\cap M^-}K dA_s=-2\pi\chi(M_\eps\cap M^-)+\int_{\partial(M_\eps\cap M^-)}k_g d\sigma.
 $$
  Hence letting $\eps$ go to zero, to show formula \r{eq:gb1} it is sufficient to prove that the contribution at boundaries offset each other, i.e.,
\begin{equation}\label{eq:offset}
  \lim_{\eps\to0}\left(\int_{M_\eps\cap M^+}k_gd\sigma-\int_{M_\eps\cap M^-}k_g d\sigma\right)=0.
  \end{equation}
   This is the key point  where the absence of tangency point is used. Indeed, under this assumption, the set $\partial M_\eps$ is shown to be smooth and, moreover, the contributions at the boundaries have the same order in $\eps$. In next section we see that this symmetry between the two sides of the singular set is lost in a neighborhood of a tangency point.

  \smallskip
  
 Theorem~\ref{th:gb1} was generalized in \cite{high-order}
  to surfaces with boundary. In this case,  the authors define admissible domains as   open bounded connected domains $U\subset M$ whose boundary is the finite union of $\con^2$-smooth admissible curves and study the convergence of $\int_{U\cap M_\eps} K dA_s$ as $\eps$ goes to zero. The generalized Gauss--Bonnet formula, proved   through the above techniques, takes account of the boundary contributions and some other terms due to the intersections of $\partial U$ with the singular set. When the intersection between $\partial U$ and $\Zz$ is $\con^2$-smooth, the formula simplifies to the following one.
  
  \bt[\cite{high-order}]\label{th:gb2} If there are no tangency points and $U$ is an admissible domain such that $\partial U$ is piecewise $\con^2$-smooth and $\partial U\cap \Zz$ is  $\con^2$-smooth, then
  $$
  \int_U K dA_s+\int_{\partial U}k_g d\sigma=2\pi(\chi(U^+)-\chi(U^-)).
  $$
 \et
   In the previous statement we have
 \begin{eqnarray*}
 \int_U K dA_s&=&\lim_{\eps\to0}\int_{U\cap M_\eps} K dA_s,\\
\int_{\partial U} K dA_s&=&\lim_{\eps\to0}\left(\int_{\partial U\cap \partial U_\eps^+}k_g dA_s-\int_{\partial U\cap \partial U_\eps^-}k_g dA_s\right),
 \end{eqnarray*}
where  $U$ has the orientation induced as a domain of $M$ and $\partial U$ is oriented as  boundary of $U$;  $U_\eps^\pm$ denotes the set $U\cap M_\eps\cap M^\pm$ and $k_g$ denotes the Gaussian curvature and $d\sigma$ the arclength parameter\footnote{Each $\con^2$-smooth piece of $\partial U$ is an admissible curve parameterized by arclength}.

  \subsubsection{ARSs with tangency points}\label{sec:gbtang}
  
  Let us present a further generalization of Theorem~\ref{th:gb1} allowing the presence of tangency point.

 In this case, two main issues are to be considered when studying the convergence \r{eq:lim}. First, around a tangency point the boundary of the domain of integration is not smooth but has two corners (one on each side of $\Zz$, see Figure~\ref{fig:astang}). Second, the Gaussian curvature diverges in a more complicated way than at Grushin points. Indeed, while at Grushin points $K$ diverges to $-\infty$, at a tangency point $K$ diverges at $+\infty$ on some directions, and it diverges at $-\infty$ at some other ones (see Figures~\ref{fig-KRESTACCIA}, \ref{fig:div}). 
Together with the interaction between different orders in the asymptotic expansion of the almost-Riemannian distance,  these remarks possibly explain  why the limit \r{eq:lim} is not  shown to converge. Indeed,  the compensation between boundary terms \r{eq:offset} appears not to happen around tangency points. 
This has been supported by numerical simulations in \cite[Section 5.1]{euler}
for the ARS on $\R^2$ having $(1,0)$, $(y-x^2)$ as orthonormal frame. For this example, the limit \r{eq:offset} appears to diverge  near $(0,0)$ as $c/\eps$, $c$ is a positive constant.


\begin{figure}[h!]
\begin{center}
\input{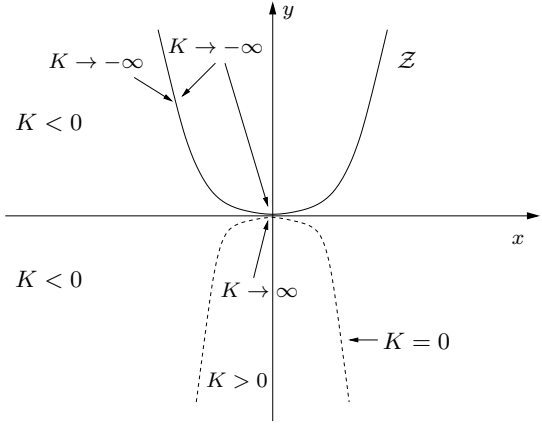}

\caption{Behavior of $K$ for the ARS on $\R^2$ whose orthonormal frame is $(1,0)$, $(0, y-x^2)$:  the plane $\R^2$ with the singular set (solid  line) and the set of points at which $K=0$ (dashed line). See also Figure~\ref{fig-KRESTACCIA} for the graph of $K$}\label{fig:div}
\end{center}
\end{figure}

To overcome the problem, a new notion of integrability of the curvature has been provided in \cite{euler}. The idea is to integrate the curvature not on the whole $M_\eps$, but on a subset of $M_\eps$ depending on two other parameters $\delta_1, \delta_2$. This set is built by taking away from $M_\eps$ a ``rectangular" box for each tangency point, where $\delta_1,\delta_2$ are the dimensions of the box see Figure~\ref{boxbr}.

\begin{figure}[h!]
\begin{center}
\input{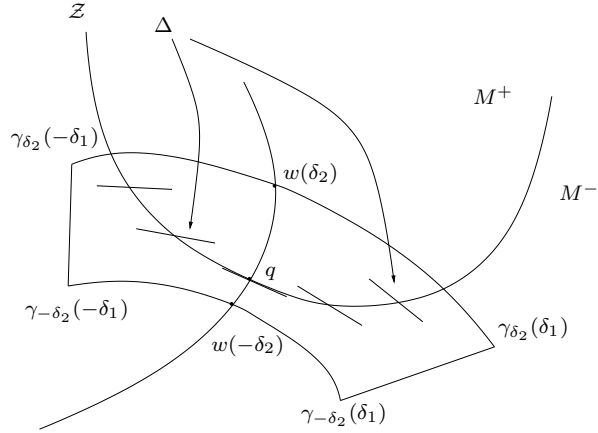}
\caption{Let $q$ be a  tangency point. Take any smooth parameterized curve 
$c:]-1,1[\to M$
such that $c(0)=q$ and $\dot c(s)\notin \bD(c(s))$. For each $s\in (-1,1)$, denote by
  $\g_s$ the  geodesic (parameterized by arclength) 
 such that $\g_s(0)=c(s)$ and $d(\g_s(t),w(]-1,1[)= |t|$ for each $t$ sufficiently small. For $\delta_1,\delta_2$ sufficiently small, the rectangle $B^q_{\delta_1,\delta_2}$ is the subset of $M$ 
 containing the tangency point $q$
 and having as boundary 
$\g_{\delta_2}([-\delta_1,\delta_1] )\cup
\g_{[-\delta_2,\delta_2]}(\delta_1)\cup
\g_{-\delta_2}([-\delta_1,\delta_1])\cup
\g_{[-\delta_2,\delta_2]}(-\delta_1)
$
 }\label{boxbr}
\end{center}
\end{figure}
 Let $M_{\eps,\delta_1,\delta_2}=M_\eps \setminus  \bigcup_{q\in {\cal T} }  B^q_{\delta_1,\delta_2}$.
  Then $K$ is said to be \emph{ 3-scale integrable with respect to the ARS} if the limit 
 \bqn\label{eq:trelimit}
\lim_{\delta_1\to0}\lim_{\delta_2\to0}\lim_{\eps\to0}\int_{M_{\eps,\delta_1,\delta_2}} K dA_s
\eqn
 exists. In this case such limit is denoted by $\treint_M K dA_s$.

\smallskip

Note that when there are no tangency points, the limit defined in~\r{eq:trelimit} clearly coincides with the one in \r{eq:lim}. Besides, the order in which the limits are taken in \r{eq:trelimit} is important. Indeed, if the order is permuted, then the result given
in Theorem~\ref{th:gb3} does not hold. 
 
 To obtain a Gauss--Bonnet formula when tangency points are present, we need the notion of contribution at tangency points. 
\bdeff\label{deftp}
Let $q$ be a tangency point. Orient $\Zz$ as the boundary of $M^+$ (see Figure~\ref{tauu}). We define the \emph{contribution at } $q$ as $\tau_q=1$, respectively $\tau_q=-1$, if the distribution is rotating counterclockwise, respectively clockwise, along $\Zz$ at $q$, see Figure~\ref{tauu}.
\edeff

\begin{figure}[h!]
\begin{center}
\input{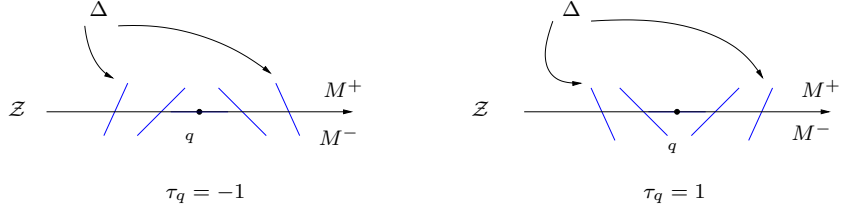}
\caption{Tangency points with opposite contributions}\label{tauu}
\end{center}
\end{figure}

The following result generalizes Theorem~\ref{th:gb1} to ARSs with tangency points.

\bt[\cite{euler}]\label{th:gb3}
The Gaussian curvature is 3-scale
integrable and
\begin{equation}\label{eq:gb3}
\treint_M K dA_s=2\pi\left(\chi(M^+)-\chi(M^-)+\sum_{q\in{\cal T}}\tau_q\right),
\end{equation}
where ${\cal T}$ is the set of tangency points of ${\cal S}$.
\et

Notice that the construction of $M_{\eps,\delta_1,\delta_2}$ depends on the choice of a manifold transversal to $\Zz$ at each tangency point and on its parameterization. A canonical choice of this manifold
has been provided in \cite[Theorem 2]{kresta}, implying that $M_{\eps,\delta_1,\delta_2}$ can be constructed in a intrinsic way. This is also suggested by the fact that  $\treint_M K dA_s$ equals a quantity that is intrinsically associated with the structure, see formula \r{eq:gb3}.

%
%

\smallskip

\noindent{\bf Proof of Theorem~\ref{th:gb3} (sketch).}
Fix $\delta_1$ and $\delta_2$ in such a way that  the rectangles 
$B^q_{\delta_1,\delta_2}$ are pairwise disjoint and
$\Zz\cap\partial B^q_{\delta_1,\delta_2}\subset [-\delta_1,\delta_1]\times\{\delta_2\}$, for every $q\in{\cal T}$.
By construction, 
$\partial B^q_{\delta_1,\delta_2}$ is admissible and has finite length for every $q\in {\cal T}$. Hence we can take  $M_\eps\setminus \bigcup_{q\in {\cal T}} B^q_{\delta_1,\delta_2}$
 as $U$ in Theorem~\ref{th:gb2}.
As a consequence, we have
 \bqn
\lim_{\eps\to0}\int_{M_{\eps,\delta_1,\delta_2}} K dA_s+\sum_{q\in{\cal T}}\int_{\partial B^q_{\delta_1,\delta_2}}k_g d\sigma_s&=&2\pi(\chi(M^+\setminus \bigcup_{q\in{\cal T}} B^q_{\delta_1,\delta_2})-  \chi( M^-\setminus \bigcup_{q\in{\cal T}} B^q_{\delta_1,\delta_2} ))\nn\\
&=&2\pi(\chi(M^+)-\chi(M^-)).
\eqn
If we prove that, for a fixed $q\in{\cal T}$, 
\begin{equation}\label{eq:secpoi}
\lim_{\delta_1\to 0}\lim_{\delta_2\to 0}\int_{\partial B^q_{\delta_1,\delta_2}}k_g d\sigma_s=-2\pi \tau_{q}, 
\end{equation}
then we directly obtain \r{eq:gb3}. To deduce \r{eq:secpoi} consider the local representation $(F3)$ and assume that  $\{(x,y)\mid y-x^2\psi(x)< 0\}\subset M^+$, the proof for the case $\{(x,y)\mid y-x^2\psi(x)< 0\}\subset M^-$ being analogous. On one hand, one can check that $\tau_{q}=1$. On the other hand, the geodesic curvature along $[-\delta_1,\delta_1]\times\{\delta_2\}$ and along $[-\delta_1,\delta_1]\times\{-\delta_2\}$ is zero, the two segments being 
the support of geodesics. Hence
$$
\int_{\partial B^q_{\delta_1,\delta_2}}k_g d\sigma_s=\int_{\{\delta_1\}\times[-\delta_2,\delta_2]}k_g d\sigma_s+
\int_{\{-\delta_1\}\times[-\delta_2,\delta_2]}k_g d\sigma_s+ \sum_{j=1}^4\alpha_j
$$ where the last term is the sum of the values of the angles of the box and is equal to $-2\pi$. Indeed, because of the diagonal form of the metric with respect to the chosen coordinates, each angle has value $-\frac \pi 2$. The first two terms are well defined and tend  to zero when $\delta_2$ tends to zero. Hence 
$$
\lim_{\delta_1\to 0}\lim_{\delta_2\to 0}\int_{\partial B^q_{\delta_1,\delta_2}}k_g d\sigma_s=-2\pi=-2\pi \tau_{q}.
$$\hfill$\blacksquare$

\subsection{A topological classification of ARSs}\label{sec:topol}

In this section we present a result that provides a relation among the topology of the almost-Riemannian surface and the Euler number  of the vector bundle associated with the structure.

Let us recall the notion of Euler number.
Given  an oriented vector bundle of rank $2$ over a compact connected oriented surface $M$,  the Euler number of $E$, denoted by $\E(E)$, is the self-intersection number of $M$ in $E$, where $M$ is identified with the zero section. To compute $\E(E)$, consider a smooth section $\sigma:M\rightarrow E$ transverse to the zero section. Then, by definition,
$$
\E(E)=\sum_{p\mid \sigma(p)=0}i(p,\sigma),
$$
where $i(p,\sigma)=1$, respectively $-1$, if $d_p\sigma:T_pM\rightarrow T_{\sigma(p)}E$ preserves, respectively reverses, the orientation. Notice that if we reverse the orientation on $M$ or on $E$  then $\E(E)$ changes sign. Hence, the Euler number of an orientable vector bundle $E$ is defined up to a sign, depending on the orientations of both $E$ and $M$. Since reversing the orientation on $M$ also reverses the orientation of $TM$, the Euler number of $TM$ is defined unambiguously and is equal to $\chi(M)$, the Euler characteristic of $M$. We refer the reader to \cite{hirsch} for a  more detailed discussion of this subject.
%
%
%
%

\bt[\cite{euler}]\label{th:topol}
Under the assumptions of Remark~\ref{assunzioni}, there holds
\begin{equation}\label{eq:topol}
\chi(M^+)-\chi(M^-)+\sum_{q\in{\cal T}}\tau_q=\E(E).
\end{equation}
\et


 The strategy to prove \r{eq:topol} is based on the construction of a section $\sigma$ of $E$ having only isolated zeros $\{p_1,\dots,p_m\}$ and such that 
 $$
 \sum_{i=1}^mi(p_i, \sigma)=\chi(M^+)-\chi(M^-)+\sum_{q\in{\cal T}}\tau_q.
 $$
 To this aim, the key point is to construct  $\sigma$ in a tubular neighborhood of $\Zz$ in such a way that it vanishes only at tangency points (see \cite[Lemma 1]{euler}) and at each $q\in{\cal T}$ the relation $i(q,\sigma)=\tau_q$ holds (see \cite[Lemma 2]{euler}). 
 Once this is done, it is sufficient to extend $\sigma$ in a smooth way to the whole surface.  By a transversality argument, this can be done by introducing only a finite number of isolated zeros whose index sum can be calculated using Hopf's Index Formula.

\smallskip

Theorem~\ref{th:topol} has several implications. First it classifies the ARS with respect to the associated vector bundle.
Indeed, the Euler number represents the only topological invariant of an oriented rank-2 vector bundle over a compact oriented surface, i.e., it  identifies the vector bundle. In particular, as a direct consequence we get that an ARS is trivializable, i.e., $E$ is isomorphic to the trivial bundle, if and only if 
$$
\chi(M^+)-\chi(M^-)+\sum_{q\in{\cal T}}\tau_q=0.
$$ This generalizes and provides the converse result of  \cite[Lemma 5]{ABS} stating that if tangency point are absent, i.e., ${\cal T}=\emptyset$, and the structure is trivializable then $\chi(M^+)-\chi(M^-)=0$. An alternative proof of the fact that the latter condition is sufficient for the structure to be trivializable can be found in \cite{proc-cdc}.

Moreover, Theorem~\ref{th:topol} allows to rewrite the Gauss--Bonnet formula of Theorem~\ref{th:gb3} as follows.
\bc\label{th:totale}
For any totally oriented ARS ${\cal S}$ on a compact surface $M$, the Gaussian curvature is $3$-scale integrable and
\begin{equation}\label{eq:totale}
\treint_M K dA_s=2\pi\e(E).
\end{equation}
\ec
 Corollary~\ref{th:totale} encloses previous Gauss--Bonnet formulas:  for Riemannian structures, where $\f$ is an isomorphism and $|\e(E)|=\e(TM)=\chi(M)$ (formula \r{eq:gbr}) and for ARSs without tangency points, where ${\cal T}=\emptyset$ (formula \r{eq:gb1}).  
 
 When considering Riemannian structures formula \r{eq:totale} says that the total curvature is zero if and only $\chi(M)=0$, that is, $M$ is diffeomorphic to the torus. Instead, in the almost-Riemannian context formula \r{eq:totale} implies that if the total curvature is zero then the vector bundle $E$ is trivial. As a consequence, there exist ARSs on surfaces of positive genus having zero total curvature, see for example Section~\ref{sec-sfera}.

\subsection{Lipschitz equivalence of ARSs}\label{sec:lipeq}

This section is  devoted to the description of  how the presence of the singular set and, in particular,  of tangency points affect the distance associated with the ARS.  Namely,
we focus our attention on the problem of Lipschitz equivalence among different almost-Riemannian distances. 

A  \emph{Lipschitz equivalence}  is a diffeomorphism $\varphi:M_1\rightarrow M_2$ which is bi-Lipschitz as a map from the metric space $(M_1,d_1)$ to $(M_2,d_2)$, where $d_i$ is an almost-Riemannian distance on the surface $M_i$ associated with an ARS $(E_i,\f_i,\ps_i)$ on $M_i$. Recall that bi-Lipschitz means that there exists a constant $C\geq 1$ such that
$$
\frac{1}{C}\,d_2(\varphi(q),\varphi(p))\leq d_1(q,p)\leq C\, d_2(\varphi(q),\varphi(p)),\quad \forall\, q, p\in M_1.
$$

In the Riemannian
case,  all  distances on diffeomorphic   compact oriented surfaces are Lipschitz equivalent. In other words the Lipschitz classification of Riemannian distances on compact oriented surfaces coincides with the differential one. 

In the almost-Riemannian case the Lipschitz classification is finer. Clearly this is due to the presence of a singular set and mainly it is due to  how the singular set splits the surface together with the location of   tangency points with their contributions. 
\subsubsection{Graph of a totally oriented ARS}
 It turns out that all the information needed to identify the Lipschitz equivalence class
of an almost-Riemannian distance  can be encoded in a labelled graph that is naturally associated with the structure. 

The  vertices of such graph correspond to 
connected components of $M\setminus\Zz$ and the edges correspond to connected components of $\Zz$. 
The edge corresponding to a connected component ${\cal W}$ of $\Zz$  joins the two vertices corresponding to the connected components of $M\setminus\Zz$ adjacent to ${\cal W}$. 
Every vertex $v$ corresponding to a component $M_v$ is labelled with a pair of integers $(\sign(v),\chi(v))$, where $\sign(v)$ takes into account of the orientation of $M_v$ ($\sign(v)=\pm 1$ if $M_v\subset M^\pm$) and $\chi(v)$ is the Euler characteristic of $M_v$.
Every edge $e$ corresponding to a component ${\cal W}\subset\Zz$ is labelled with the ordered sequence  of signs (modulo cyclic permutations) given by the contributions at the tangency points belonging to ${\cal W}$ where the order is fixed by walking along ${\cal W}$ oriented as the boundary of $M^+$. 
In Figure~\ref{fig:come} we illustrate the algorythm to build the labelled graph associated with an ARS.

\begin{figure}[h!5cm]
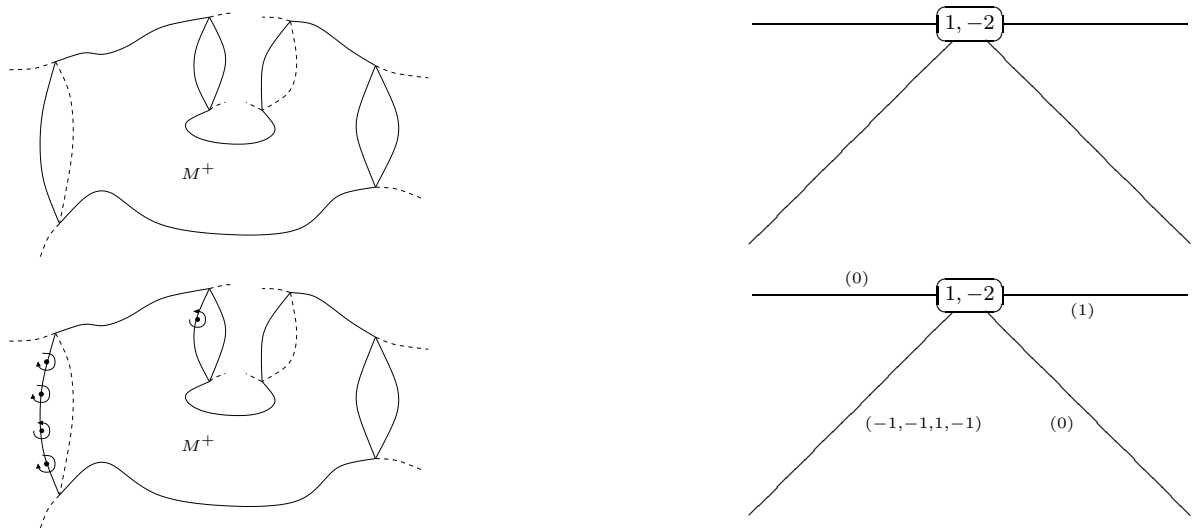

\scriptsize{
\xymatrix@=3cm@!0{*!<4cm,1.5cm>{\input{punto1.pstex_t}}&&\ver{1,-2}\ar@{-}[l]\ar@{-}[r]\ar@{-}[dr]\ar@{-}[dl]&\\
&&&\\
}
\xymatrix@=3cm@!0{*!<4cm,1.5cm>{\input{punto2.pstex_t}}&&\ver{1,-2}\ar@{-}[l]_{(0)}\ar@{-}[r]_{(1)}\ar@{-}[dr]_{(0)}\ar@{-}[dl]^{(-1,-1,1,-1)}&\\
&&&}
}
\caption{\label{fig:come}Algorythm to build the graph}
\end{figure}

An example of ARS and associated graph is portraited\footnote{Note that in Figure~\ref{fig:exg} we represent only the surface,  the singular set and the tangency points. Here we are not interested in the explicit expression of the morphism $\f$ nor on the vector bundle $E$.} in Figure~\ref{fig:exg}. In this case we consider an ARS on  the compact oriented surface of genus 2 where $\Zz$ is the union of two circles placed as in Figure~\ref{fig:exg}a and where there is only a tangency point of positive contribution on one of them. Then the graph (Figure~\ref{fig:exg}b) must have 3 vertices and 2 edges. The label on  vertices are easily computed. As concerns labels on edges there is only a tangency point (with positive contribution) on a connected component of $\Zz$, whence only one edge carries a label $(1)$, while the other edge is labelled with $(0)$. 

\brem\label{rk:or}
Adding the quantity $\sum_v\sign(v)\chi(v)$ to the sum of all the entries in the label of edges,   we obtain $\chi(M^+)-\chi(M^-)+\sum_{q\in{\cal T}}\tau_q$, which equals $\e(E)$, by Theorem~\ref{th:topol}. Also, $\sum_v\chi(v)=\chi(M)$. Hence, once the labelled graph associated with an ARS is given, one  recovers directly  the vector bundle and the surface.

Moreover, the labelled graph associated with an ARS depends on the orientation fixed on $E$. More precisely, choosing on $E$ the opposite orientation produces the following changes in the labels of the graph. On each vertex the first entry of the label changes sign. On each edge not only each entry of the tuple changes sign but also the tuple is reversed in order.\footnote{For example if the label of an edge is $(-1,1,-1,-1)$, changing the orientation on $E$ the new label becomes $(1,1,-1,1)$.}
\erem

\begin{figure}[h!]
\begin{center}
\input{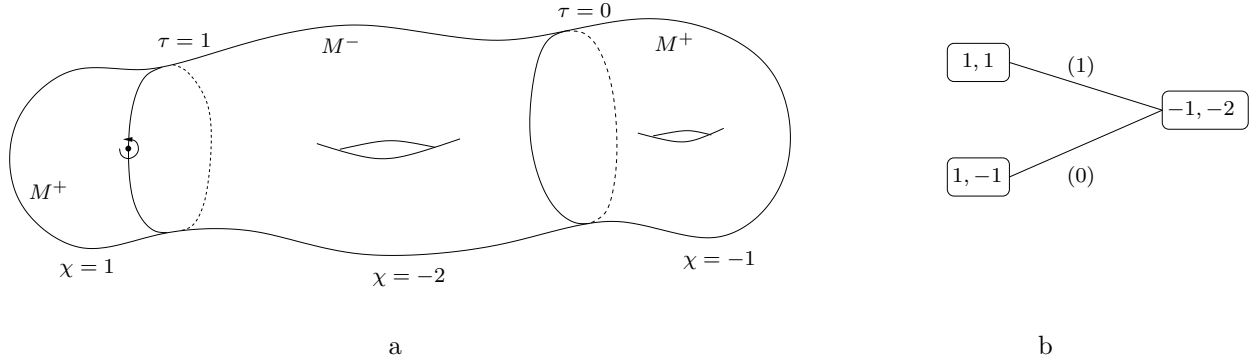}
\end{center}
\caption{Example of ARS and associated graph}\label{fig:exg}
\end{figure}

We say that two labelled graphs associated with two totally oriented ARSs are \emph{equivalent} if either they are equal or after possibly changing the orientation on one vector bundle they are equal.
This notion of equivalence is motivated by this  straightforward remark: changing the orientation of the vector bundle does not affect the almost-Riemannian distance. Moreover, when two labelled graphs are equivalent, by Remark~\ref{rk:or} the associated vector bundles are isomorphic and the underlying surfaces are diffeomorphic.

\subsubsection{A classification result}

The following result classifies totally oriented ARSs with respect to Lipschitz equivalence. 
\bt[\cite{BCGS}]\label{lip-eq}
 Two totally oriented almost-Riemannian structures defined on compact surfaces  are Lipschitz equivalent if and only if they have equivalent graphs.
\et

Theorem~\ref{lip-eq}  implies that the Lipschitz equivalence class of an almost-Riemannian distance on a surface depends only on how the singular set is embedded in the surface and how the tangency points with their contribution are located. 

In the Riemannian context\footnote{If the structure is Riemannian, the associated labelled graph consists of a unique vertex and has no edges. The label on the vertex is $(\delta,\chi(M))$, where $\delta=1$ if $\f$ preserves the orientation, $-1$ otherwise.}, the Lipschitz equivalence class of distances on a given surface is unique or, equivalently, it does not depend on the bilinear form on the tangent bundle. In a similar way, by Theorem~\ref{lip-eq} we deduce that the Lipschitz equivalence between two distances does not depend on the bilinear form $G$ defined on $\bD$ (see Section~\ref{sec:prel}) but only on the submodule $\Delta$ itself. This is highlightened by the fact that the graph itself depends only on $\Delta$. In terms of the triple $(E,\f,\ps)$ this translates to the fact that the labelled graph does not depend on the chosen Euclidean structure $\ps$, but only on  the morphism $\f$.
As a consequence, in general Lipschitz equivalence does not imply isometry. 

\smallskip

The main tool in the  proof of Theorem~\ref{lip-eq} is a local classification of ARS by Lispchitz equivalence. Indeed, using the Ball-Box Theorem (see \cite[Corollary 7.35]{bellaiche}) it is not hard to show that in a neighborhood of a given point two ARSs are Lipschitz equivalent if and only if the point is of the same type (Riemannian, Grushin, tangency) for both structures. Then, to glue the information on different neighborhoods 
one needs the topological classification of surfaces, that is, one uses the information carried by labels on vertices. 
To appreciate the role of contributions of tangency points, recall that in a neighborhood of a tangency point the asymptotic of the almost-Riemannian distance from $\Zz$ is  different from the two sides of  $\Zz$, see Section~\ref{sec:locres}. The contribution $\tau_q$ carries the information about which side is the one with the fastest rate of convergence of $M_\eps$ (see \r{eq:meps}) to $\Zz$.


\section{Conclusions}\label{sec:riemvsariem}

In Table~1 we compair Riemannian and almost-Riemannian geometry of surfaces summing up the main aspects presented in this paper.

\begin{figure}[h!]
\begin{center}
\input{resume.pstex_t}
\end{center}
\end{figure}

As one can infer from the present analysis,
the most interesting points of ARSs are tangency points. Even though some contributions have been done \cite{BCGJ,kresta}, tangency points are far to be deeply understood.  

An open question arisen in the proof of Theorem~\ref{th:gb1} is the convergence or the divergence of the integral of the geodesic curvature on the boundary of  a tubular  neighborhood  of the singular set, close to a tangency point.
To address this issue one needs to understand whether the domain of integration or the form to be integrated need to be reconsidered. 

The analysis of the Laplace--Beltrami operator in presence of tangency points has not been considered yet although in \cite{camillo} the authors conjecture some properties, based on the case with only Grushin points.

Another interesting problem is whether it is possible to associate with an ARS a canonical linear connection on $E$ compatible with the Euclidean structure. To this aim one needs to use properties of the morphism $\f$, as in general the answer is negative.
This could be a step forward towards an intrinsic notion of integration of the Gaussian curvature on the surface, as one could study the curvature on the vector bundle. Also, focusing on the vector bundle rather than on the (tangent bundle to the) surface is a starting point towards higher dimensions generalizations.

\small{
\bibliographystyle{abbrv}
\bibliography{biblio_grenoble}
}

\end{document}